\documentclass[10pt, a4paper]{amsart}

\usepackage{amsmath}
\usepackage{wasysym}
\usepackage{amsthm}
\usepackage{latexsym}
\usepackage{amsfonts}
\usepackage{mathrsfs}
\usepackage{amssymb}
\usepackage{dsfont}
\usepackage[all]{xy}
\usepackage{color}
\usepackage{cleveref}

\newtheorem{Theorems1}{Theorem}[section]

\newtheorem{Coroll1}[Theorems1]{Corollary}
\newtheorem{Lemma1}[Theorems1]{Lemma}
\newtheorem{Proposition1}[Theorems1]{Proposition}

\theoremstyle{definition}
\newtheorem{Definitions1}{Definition}[section]

\theoremstyle{plain}

\theoremstyle{remark}

\theoremstyle{plain}

\newtheorem*{Ques}{Question}

\newtheorem*{Theo}{Theorem}

\newtheorem*{Cor}{Corollary}

\theoremstyle{remark}
\newtheorem{Exam}{Example}[section]
\theoremstyle{plain}

\theoremstyle{remark}
\newtheorem{Rem}{Remark}[section]

\numberwithin{equation}{Theorems1}

\numberwithin{equation}{Theorems1}
\def\multiset#1#2{\ensuremath{\left(\kern-.3em\left(\genfrac{}{}{0pt}{}{#1}{#2}\right)\kern-.3em\right)}}

\def\Multiset#1#2{\left(\!\!\left(\!\!
    \begin{array}{c}
      n \\
      r
    \end{array}
  \!\!\right)\!\!\right)}

\newcommand{\depth}[1]{\operatorname{depth}(#1)}
\newcommand{\spec}[1]{\operatorname{Spec}(#1)}
\newcommand{\pd}[1]{\operatorname{pd}(#1)}
\newcommand{\id}[1]{\operatorname{id}(#1)}
\newcommand{\fd}[1]{\operatorname{fd}(#1)}

\newcommand{\supp}[1]{\operatorname{Supp}(#1)}

\newcommand{\rank}[1]{\operatorname{rank}(#1)}


\begin{document}

\title{Ore Localization and Minimal Injective Resolutions}
\author{Rishi Vyas}
\address{
Rishi Vyas, Wolfson College, Cambridge, CB3 9BB, United Kingdom}
\email{rv251@cam.ac.uk}
\date{}

\begin{abstract}
In this paper, we describe the structure of the localization of $\mathrm{Ext}^{i}_{R}(R/P,M)$, where $P$ is a prime ideal and $M$ is a module, at certain Ore sets $X$. We first study the situation for FBN rings, and then consider matters for a large class of Auslander-Gorenstein rings. We need to impose suitable homological regularity conditions to get results in the more general situation. The results obtained are then used to study the shape of minimal injective resolutions of modules over noetherian rings.  
\end{abstract}

\maketitle

If $R$ is a commutative noetherian ring, $S\subseteq R$ a multiplicatively closed set containing $1$, and $M$ a finitely generated module, then for any $R$-module $N$ and any integer $i$, there is an isomorphism: $$\mathrm{Ext}^{i}_{R}(M,N)_{S}\cong \mathrm{Ext}^{i}_{R_{S}}(M_{S},N_{S}).$$ In general, the question of whether there is an analogous result for noncommutative noetherian rings does not make sense; $\mathrm{Hom}_{R}(M,N)$ need not carry an $R$-module structure for arbitrary modules $M$ and $N$. However, if $M$ or $N$ is a bimodule, then such questions become valid and can be quite subtle to answer. 

In this paper, we will show that for some rings $R$ at certain Ore sets $X$, given a finitely generated module $M$ and a prime ideal $P$, there is an isomorphism $\mathrm{Ext}^{i}_{R}(R/P,M)\otimes_{R} R_{X}\cong \mathrm{Ext}^{i}_{R_{X}}(R_{X}/P_{X},M_{X})$. This is enough for some applications. Given a module $M$ over a noetherian ring $R$, what can one say about the `shape' of a minimal injective resolution of $M$? For a commutative noetherian ring $A$, this is known (see \cite{ROBS}, \cite{FFGR}, \cite{FOX1}): the minimal injective resolution of a module $M$ is determined by geometric data corresponding to $M$, i.e. the support of $M$, $\operatorname{Supp}(M)$, in $\operatorname{Spec}(A)$. Using the results in this paper, we will show that matters behave similarly for some classes of noncommutative noetherian rings.

In \textsection $1$, we gather some preliminary results and definitions. 

In \textsection $2$,  we make our first attempt to localize $\mathrm{Ext}$. We show that if $X$ is an Ore set such that essentiality is preserved under localization (such Ore sets are characterized in \cite{GJ2}), then matters behave well without additional homological hypotheses. Over an FBN ring, localization at every Ore set preserves essentiality --- this was first proved in \cite{GJ}, and we provide a proof of a slightly more general result in \cref{localizeinj}. \cref{localizeext}, in particular, proves the following:

\begin{Theo}
Let $R$ be a noetherian right FBN ring. Let $X$ be an Ore set in $R$. Let $M$ be an $R$-module. Then, given a prime ideal $P$, there is an isomorphism of $R_{X}$-modules for all $i$: $$ \mathrm{Ext}^{i}_{R}(R/P,M)\otimes_{R} R_{X}\cong \mathrm{Ext}^{i}_{R_{X}}(R_{X}/P_{X},M_{X}).$$ If, in addition, $S$ is a ring and $M$ is an $(S,R)$-bimodule, then the above isomorphism is one of bimodules. 
\end{Theo}

In \cite[Lemma 3.2]{KB}, this was proved, with a mild homological hypothesis, for bimodules. We reprove this result without the homological restrictions on $R$.

Recall that over a right noetherian ring, indecomposable injective modules come in two distinct flavours, tame and wild (see \textsection 1.1 for details). Fully Bounded Noetherian (FBN) rings are precisely the noetherian rings for which all indecomposable injectives are tame, and this fact allows us to work directly with injective resolutions to prove localization results. Injective modules over rings which are not FBN, however, can be far more difficult to understand, and in this case our direct approach fails miserably. 

When working with noetherian rings which are not FBN, we will need to impose strong homological restrictions and restrict the Ore sets we are localizing at to those which arise from cliques of localizable prime ideals. The key idea here is to approximate $\mathrm{Ext}^{i}_{R}(R/P,M)$ by a spectral sequence, which is then easier to localize. In order to do this, we need to understand the structure of $\mathrm{Ext}^{i}_{R}(R/P,R)$ as an $(R,R)$-bimodule. In \cite[Theorem 5.5]{BL} Brown and Levasseur studied this very object for a prime ideal $P$ of height $t$ in the universal enveloping algebra $U$ of a soluble Lie algebra $\mathfrak{g}$, and proved that in this case $\mathrm{Ext}_{U}^{t}(U/P,U)$ is isomorphic to a right ideal of $U/P$ as a right module and is isomorphic to a left ideal of $U/\tau^{*}(P)$ as a left module, where $\tau^{*}$ is the winding automorphism induced by a certain character $\tau$ of $\mathfrak{g}$. 

In \textsection{3} of this paper, we aim to extend the above result to a wider class of rings. Our general result will have the following form: for a suitably nice ring $R$, prime ideal $P$, and $i$, there exists another prime ideal $Q$ such that $\mathrm{Ext}^{i}_{R}(R/P,R)$ is annihilated as a left module by $Q$, and is torsionfree both as a right $R/P$-module and a left $R/Q$-module.

We will describe two ways of doing this. The first is a somewhat elementary argument that relies heavily on Ore localization; it only applies to a reasonably restrictive class of rings - Auslander-Gorenstein, grade-symmetric, weakly-bifinite rings which satisfy the second layer condition. For the second, the author is indebted to James Zhang, and his observation that the theory of dualizing complexes can be used to prove and extend the result in many situations. Of particular interest are \cref{bimodstruc} and \cref{bimodstruc2} (note the paragraph after \cref{dual1}).

In \textsection 4, we use the results of \textsection 3 along with a spectral sequence (\cref{iss1}) to localize $\mathrm{Ext}^{i}_{R}(R/P,M)$, where $M$ is a module of finite flat dimension. The following is a special case of \cref{locmod}:

\begin{Theo}
Let $R$ satisfy one of the following conditions:
\begin{itemize}
\item[a.)] $R$ is an Auslander-Gorenstein, grade-symmetric, weakly bifinite ring which satisfies the second layer condition.
\item[b.)] $R$ is a $k$-algebra, where $k$ is a field, with a noetherian, connected filtration such that $gr(R)$ is a commutative AS-Gorenstein ring.
\item[c.)] $R$ is an Auslander-Gorenstein, complete local $k$-algebra with maximal ideal $J$ such that $R/J$ is finite dimensional over $k$.
\end{itemize}
Let $P$ be a prime ideal with localizable clique $C$ with corresponding Ore set $X$. Then, given a module of finite flat dimension $M$, there is an isomorphism of $R/P$-modules: $$\mathrm{Ext}^{i}_{R}(R/P,M)\otimes_{R} R_{X}\cong \mathrm{Ext}^{i}_{R_{X}}(R_{X}/P_{X},M_{X}).$$
\end{Theo}

There are several results describing the structure of a minimal injective resolution of a ring over itself, in terms of multiplicities of indecomposable injectives (see \cite{JY} for references to the literature and state of the art technology) and with regards to purity under various dimension functions (\cite{ASZ}, \cite{ASZ2}). In \textsection 5, the results in this paper are used to study the `shape' of a minimal injective resolution of a module over certain noetherian rings. \cref{resin1}, \cref{resin2}, and \cref{finmult} combine to prove the following:

\begin{Cor}
Let $k$ be an uncountable field, and let $A$ be a $k$-algebra  of injective dimension $d$ which satisfies one of the following conditions:
\begin{itemize}
\item[a.)] $A$ has a noetherian, connected filtration such that $gr(A)$ is a commutative AS-Gorenstein ring.
\item[b.)] $A$ is a Auslander-Gorenstein, complete local ring with maximal ideal $J$ such that $A/J$ has finite dimension over $k$.
\end{itemize}
Let $M$ be a finitely generated module of finite injective dimension with minimal injective resolution $\mathrm{\mathbf{I}}:= 0\rightarrow I^{0}\rightarrow I^{1}\rightarrow ...$, and let $C$ be a localizable clique with corresponding Ore set $X$. Given a prime ideal $P$, let $E_{P}$ be the unique tame indecomposable injective module with assassinator $P$. Then, $M_{X}\neq 0$ if and only if there exists a $P\in C$ such that $E_{P}$ appears as a summand of some term in $\mathrm{\mathbf{I}}$. In this case, given $\depth{M_{X}}\leq i\leq j_{R}(R/P)$, there exists $Q\in C$ such that $E_{Q}$ appears as a summand of $I^{i}$, and the multiplicity of $E_{Q}$ is finite.
\end{Cor}

We assume that the reader is familiar with Auslander-Gorenstein rings and dualizing complexes; at the beginning of \textsection 1.4, we give some references for details and definitions about these objects.  We also assume familiarity with the theory of rings which satisfy the second layer condition and the theory of Ore  localization at a prime ideal in a noncommutative ring. A good reference for this material is \cite{GW}.

We will always assume that our rings have $1$, and one-sided modules will always be right modules, unless we explicitly mention otherwise. We will always compute $\mathrm{Ext}$ in the category of right modules. All the rings in this paper will be assumed to be right noetherian; however, we will always give precise hypotheses in our results, for the sake of clarity. We will use the notation $R$ to denote an arbitrary ring satisfying some given hypotheses; however in \textsection $3.2$ we will use $A$ instead, to remind ourselves that we are working with algebras over a fixed commutative ring.

\section{Some Preliminaries and Notation}

\subsection*{Injective Modules over Noetherian Rings}

Let $R$ be a right noetherian ring. Every injective module is a direct sum of indecomposable injective modules (\cite[Theorem 3.48]{LAM}). It is well known that every indecomposable injective module $E$ has a unique associated prime, which is then known as the assassinator of $E$ (\cite[Lemma 5.26]{GW}). Given an indecomposable injective module $E$ with assassinator $P$, it is a fact that $\mathrm{ann}_{E}(P)$ is either torsionfree or torsion as an $R/P$-module (\cite[Proposition 7.10]{GW}). In the first case, $E$ is called \textit{tame}; otherwise, we call it \textit{wild}. It can be shown that given a prime $P$, there exists a unique tame indecomposable injective with associated prime $P$: we denote this distinguished injective module by $E_{P}$. In this case, $\mathrm{ann}_{E_{P}}(P)$ is easy to describe: it is the unique simple module over $\mathcal{Q}(R/P)$, the Goldie quotient ring of $R/P$. 

We will use $\id{M}$, $\pd{M}$, and $\fd{M}$ to denote the injective, projective, and flat dimension of a module $M$, respectively.

\subsection*{Noncommutative Localization}

Over a noncommutative ring, it is not always possible to localize at a prime ideal. Given a noetherian ring $R$, there is a relation on $\spec{R}$, known as a \textit{link} and denoted by $\rightsquigarrow$ (see \cite[Chapter 12]{GW}). The connected components of $\spec{R}$ with regards to this relation are known as \textit{cliques}. To every clique in a noetherian ring we can associate a multiplicatively closed set; in particular, there is the notion of a \textit{localizable} clique. Recall that if $M$ is an $R$-module, the \textit{socle} series of $M$ is the ascending sequence of subgroups defined inductively as $\operatorname{soc}^{n}(M)/\operatorname{soc}^{n-1}(M) = \operatorname{soc}(M/\operatorname{soc}^{n-1}(M)).$

 \begin{Definitions1} \label{localizableclique}
 Let $\mathrm{C}$ be a clique of prime ideals in a noetherian ring $R$. We say that $\mathrm{C}$ is localizable if $X$, the intersection of elements regular modulo $Q$ for every prime $Q\in \mathrm{C}$, is a right and left Ore set, and the localization $R_{X}$ satisfies the following properties:
\begin{itemize}
 \item[a.)] Prime ideals of the form $Q_{X}$, as $Q$ varies across $\mathrm{C}$, are the only right and left primitive ideals in $R_{X}$.
 \item[b.)] The quotient ring $R_{X}/Q_{X}$ is simple artinian for all $Q\in C$.
\end{itemize}
In addition, if for all $Q\in C$, $E(R_{X}/Q_{X})$, the injective hull of $R_{X}/Q_{X}$ as an $R_{X}$-module, is the union of its socle series, we say that $C$ is classically localizable.
\end{Definitions1}

As the only left and right primitive ideals of $R_{X}$ are those of the form $Q_{X}$ for $Q\in C$, we have that the left and right simple $R_{X}$-modules are parametrized by the prime ideals in the clique. For the purposes of this paper, we will use the notation $S^{P}$ to refer to the unique simple right $R_{X}$-module with annihilator $P_{X}$, and $^{P}S$ for the unique simple left $R_{X}$-module with the same property.

The ability to localize at a clique is strongly connected to the \textit{second layer condition}. For more details about this condition, and other aspects of localization at a prime ideal in a noncommutative ring, the reader may refer to \cite[Chapters 12, 14]{GW}.

\subsection*{Elementary Homological Identities} \label{elemhomoiden}

We state some elementary homological identities. If $M$ is a right $R$-module, and $\phi$ is an automorphism of $R$, then we denote by $M^{\phi}$ the $R$-module with the same underlying abelian group as $M$, and action $m\cdot r = m\phi(r)$. We similarly have $^{\phi}M$ for left modules. Note that  $(M^{\phi})^{\rho}\cong M^{\phi \rho}$ while $^{\phi}(^{\rho}M)\cong {^{\rho \phi}M}$. Also note that  $^{\phi}R^{1}\cong {^{1}R^{\phi^{-1}}}$, via the map $\phi$.  As mentioned earlier, $\mathrm{Ext}^{i}_{R}(M,N)$ in general does not carry a module structure, unless $M$ or $N$ happens to be a bimodule. If $M$ is an $(S,R)$-bimodule, then $\mathrm{Ext}^{i}_{R}(M,N)$ is a right $S$-module, while if $N$ is an $(S,R)$-bimodule, then $\mathrm{Ext}^{i}_{R}(M,N)$ is a left $S$-module. The following identities are easily verified by direct computation.

\begin{Lemma1} \label{iden1}
Let $M$ and $N$ be right $R$-modules. Let $\phi$ be an automorphism of $S$. Then, we have the following identities:
\begin{itemize}
 \item[a.)] Let $N$ be an $(S,R)$-bimodule. Then, $\mathrm{Ext}^{i}_{R}(M,^{\phi}N)\cong \mathrm{^{\phi}Ext}^{i}_{R}(M,N)$.
 \item[b.)] Let $M$ be an $(S,R)$-bimodule. Then, $\mathrm{Ext}^{i}_{R}(^{\phi}M,N)\cong \mathrm{Ext}^{i}_{R}(M,N)^{\phi}$.
\end{itemize}
\end{Lemma1}

\subsection*{The Ischebeck Spectral Sequences}

The following spectral sequences were originally discovered by Ischebeck in \cite{Ish}. They are also described by Kr\"{a}hmer in \cite{U}. Proofs of the following two lemmas can be found in \cite{Ish} or in \cite{U}.

\begin{Lemma1}\label{iss1}
Let $R$ be a right noetherian ring. Let $S$ and $T$ be rings. Let $_{S}M_{R}$ and $_{T}N_{R}$  be bimodules, with $_{S}M_{R}$ finitely generated as a right $R$-module. Suppose either $M$ has finite projective dimension as a right $R$-module, or $N$ has finite flat dimension as a right $R$-module. Then, there exists a convergent spectral sequence of $(T,S)$-bimodules:
$$E^{2}_{p,q}\cong \mathrm{Tor}_{p}^{R}(N, \mathrm{Ext}_{R}^{q}(M,R))\Rightarrow \mathrm{Ext}_{R}^{q-p}(M,N).$$
\end{Lemma1}

\begin{Lemma1} \label{iss2}
Let $R$ be a right noetherian ring. Let $S$ and $T$ be rings. Let $_{S}M_{R}$ and $_{R}N_{T}$ be bimodules, with $_{S}M_{R}$ finitely generated as a right $R$-module. Suppose either $M$ has finite projective dimension as a right $R$-module, or $N$ has finite injective dimension as a left $R$-module. Then, there exists a convergent spectral sequence of $(S,T)$-bimodules:
$$E_{2}^{p,q}\cong \mathrm{Ext}^{p}_{R}(\mathrm{Ext}^{q}_{R}(M,R),N)\Rightarrow \mathrm{Tor}_{q-p}^{R}(M,N).$$  
\end{Lemma1}

We make an observation: applying an exact functor does not damage the convergence of a convergent spectral sequence. To be precise, we have the following lemma, which follows immediately from the definitions of an exact functor, a spectral sequence, and the convergence of a spectral sequence. 

\begin{Lemma1} \label{flatspecseq}
Let $\mathcal{A}$ and $\mathcal{B}$ be abelian categories, and let $F:\mathcal{A}\rightarrow \mathcal{B}$ be an exact functor. Let $\mathrm{E_{p,q}^{r}}$ be a bounded spectral sequence converging to a graded object $H^{n}$. Then, $F(\mathrm{E_{p,q}^{r}})$ is a bounded spectral sequence converging to $F(H^{n})$.
\end{Lemma1}

 Thus, in the context of \cref{iss1}, if $M$ is an $(S,R)$-bimodule, and $L$ is a flat right $S$-module, we end up with a convergent spectral sequence $$E_{2}^{p,q}\cong \mathrm{Tor}_{p}^{R}(N, \mathrm{Ext}_{R}^{q}(M,R))\otimes_{S} L\Rightarrow \mathrm{Ext}_{R}^{q-p}(M,N)\otimes_{S} L.$$

Over rings of finite injective dimension, modules of finite projective dimension often coincide with modules of finite injective dimension. This has been noticed in \cite[The proof of Lemma 5.12]{LEV} and \cite[The comment after Lemma 2.1]{QZ1}. Using the two spectral sequences in \cref{iss1} and \cref{iss2}, one can prove a more precise statement. The following must be well known, but we have not seen it in the literature. So we sketch a proof:

\begin{Proposition1} \label{equiinjf}
Let $R$ be a right noetherian ring. Then, the following are equivalent:
\begin{itemize}
\item [a.)] $R$ has finite injective dimension as a right module.
\item [b.)] Every flat right module has finite injective dimension.
\item [c.)] Every injective left  module has finite flat dimension.
\end{itemize}
\begin{proof}
The fact that $a.)\implies b.)$ and $a.)\implies c.)$ are consequences of the spectral sequences of \cref{iss1} and  \cref{iss2} collapsing. Indeed, suppose $a.)$ holds. Given a flat right module $N$, the spectral sequence in \cref{iss1} collapses to give an isomorphism $$N\otimes_{R} \mathrm{Ext}^{i}_{R}(M,R)\cong \mathrm{Ext}^{i}_{R}(M,N)$$ for any finitely generated right module $M$. Since $R$ has finite injective dimension, it follows that $N$ has finite injective dimension; indeed, it follows that the injective dimension of $N$ is bounded by the injective dimension of $R$. Thus $b.)$ holds. 

Similar considerations with the spectral sequence in \cref{iss2} tell us that $a.)\implies c.)$.

It is immediate that $b.)\implies a.)$. Suppose $c.)$ holds. If $E$ is any left injective module and $M$ is any right $R$-module, \cref{iss2} gives an isomorphism: $$\mathrm{Hom}_{R}(\mathrm{Ext}^{i}_{R}(M,R),E)\cong \mathrm{Tor}_{i}^{R}(M,E).$$ Let $E'$ be an injective cogenerator in the category of left $R$-modules. Since every injective left module has finite flat dimension, there exists a $k$ such that $\mathrm{Hom}_{R}(\mathrm{Ext}^{i}_{R}(M,R),E')=0$ for all $i\geq k$. This forces $\mathrm{Ext}^{i}_{R}(M,R)=0$ for all $i\gneq k$. Thus $R$ has finite injective dimension as a right module. Therefore $c.)\implies a.)$.
\end{proof}
\end{Proposition1}

\begin{Rem}
In particular, if $R$ is noetherian and has finite injective dimension on both sides, then a module has finite flat dimension if and only if it has finite injective dimension. We will use this fact without mention for the rest of this paper, choosing one or the other depending on what is more suitable in any given context.
\end{Rem}

\subsection*{A Localization Result}

The next result is well-known and occurs in various forms throughout the literature. A proof can be found at \cite[Proposition 1.6]{BL}.

\begin{Lemma1} \label{flatten}
Let $R$ be a right noetherian ring, $S$ a ring, and let $f:R\rightarrow S$ be a ring homomorphism such that $S$ is a flat left and right $R$-module. Let $M$ be a finitely generated $R$-module, and let $B$ be an $(R,R)$-bimodule such that $S\otimes_{R} B$ is an $(S,S)$-bimodule. Then, for all $i \geq 0$, we have that $$S\otimes_{R} \mathrm{Ext}^{i}_{R}(M,B)\cong \mathrm{Ext}_{S}^{i}(M\otimes_{R} S, S\otimes_{R} B).$$ The above isomorphism is one of left $S$-modules. If $T$ is a ring and $M$ is a $(T,R)$-bimodule, then the above isomorphism is one of  $(S,T)$-bimodules.
\end{Lemma1}

The proof of the above result is not difficult: essentially, localizing $\mathrm{Ext}$ when the second variable carries a bimodule structure is straightforward because we can take a free resolution of the first module and compute. Matters are much more difficult when we want to localize at the first variable, as we then ideally need to work with an injective resolution of the second module. The next section shows how this can be done when our resolution is of an especially nice form. In most cases, however, this direct computational approach is impossible, and we need to find other methods, such as those in \textsection $3$.

\subsection*{Auslander-Gorenstein Rings, Dualizing Complexes, and Purity}

 We will assume that the reader is familiar with the theory of Auslander-Gorenstein rings and its generalization, the theory of Auslander dualizing complexes. In this paper, we will often say that $R$ is a dualizing complex over $(A,B)$ to denote that it is a dualizing complex between $A$ and $B$. For definitions and details about these objects, we refer the reader to \cite{ASZ}, \cite{QZ2}, \cite{QZ1}, \cite{JY2}, and \cite{JY}. We use \cite{WIEB} as a reference for homological algebra. 

Over a right noetherian ring, there is an integer valued function $j_{R}$ on the category of finitely generated modules, known as the \textit{grade}. $$j_{R}(M) := \mathrm{inf}\{i|\mathrm{Ext}^{i}_{R}(M,R)\neq 0\}.$$

Over an Auslander-Gorenstein ring $R$ $-j_{R}$ is a dimension function, in the sense of \cite[Definition 6.8.4]{MR}. 

An $R$-module $M$ is said to be \textit{holonomic} if $j_{R}(M)$ is equal to the injective dimension of $R$ over itself. It is a fact that every holonomic module over an Auslander-Gorenstein ring is artinian (\cite[Theorem 3.4.2.3(2)]{z-f}). 

We say that a ring is \textit{grade-symmetric},  \cite{ASZ}, if given any bimodule $_{R}B_{R}$ finitely generated on both sides, $j_{R}(_{R}B)=j_{R}(B_{R})$. Borrowing from the theory of dualizing complexes, we say that an Auslander-Gorenstein ring $R$ is \textit{bifinite} (resp. \textit{weakly bifinite}) if $\mathrm{Ext^{i}_{R}(B,R)}$ is finitely generated on both sides for any $(R,R)$-bimodule $B$ (resp. subquotient of $R$) which is finitely generated on both sides. Of course, in the case of a dualizing complex $R$ over a pair of algebras $(A,B)$, we have functions $j_{R,A}$ or $j_{R,B}$ over the category of  finitely generated left $A$-modules and right $B$ modules, respectively.  There is a notion of a Cdim-symmetric (\cite[Definition 1.8]{JY}) dualizing complex, analogous to a grade-symmetric ring, and a weakly bifinite (\cite[Definition 1.3]{JY}) dualizing complex, from where we borrowed our terminology for rings. A ring $R$ is said to be \textit{AS-Gorenstein} if it is noetherian, has finite right and left injective dimension $d$, and given any simple right or left $R$-module $S$, $\mathrm{Ext}^{i}_{R}(S,R)=0$ for all $i\neq d$, and $\mathrm{Ext}^{d}_{R}(S,R)\neq 0$ is simple. The reader should note that for noetherian connected graded algebras over a field, there is another definition of AS-Gorenstein (\cite[p. 45]{JY2}).

Recall that a module $M$ is said to be \textit{pure} of grade $i$ if every finitely generated submodule of $M$ has grade $i$. A proof of the following lemma can be found at \cite[Proposition 3.4.2.8]{z-f}.

\begin{Lemma1} \label{extpure}
Let $R$ be an Auslander-Gorenstein ring, and let $M$ be a finitely generated module. Then, $\mathrm{Ext}^{j_{R}(M)}_{R}(M,R)$ is pure of grade $j_{R}(M)$.
\end{Lemma1}

\begin{Lemma1} \label{onlyone}
Let $R$ be an Auslander-Gorenstein, grade-symmetric, weakly bifinite ring. Let $P$ be a prime ideal in $R$. Then, $\mathrm{Ext}^{i}_{R}(R/P,R)$ is torsion as an $R/P$-module if and only if $i\neq j_{R}(R/P)$, and $\mathrm{Ext}^{j_{R}(R/P)}_{R}(R/P,R)$ is torsionfree as an $R/P$-module.
\begin{proof}
If $\mathrm{Ext}^{j_{R}(R/P)}_{R}(R/P,R)$ has nontrivial torsion as an $R/P$-module, we can find a sub-bimodule $B$ such that $B$ is a torsion $R/P$-module on the right, and is finitely generated on both sides. As we are assuming that $R$ is grade-symmetric, we have that $j_{R}(_{R}B)=j_{R}(B_{R})$.  Thus $j_{R}(B_{R})=j_{R}(R/P)$, because as a left module,  $\mathrm{Ext}^{j_{R}(R/P)}_{R}(R/P,R)$ is pure of grade $j_{R}(R/P)$.  This is a contradiction, as $j_{R}$ is a dimension function on the category of finitely generated modules, which forces $j_{R}(B_{R})\gneq j_{R}(R/P)$. Thus $\mathrm{Ext}^{j_{R}(R/P)}_{R}(R/P,R)$ is torsionfree as a right $R/P$-module. 

If $i\gneq j_{R}(R/P)$, then we can argue as in \cite[Theorem 4.2]{ASZ}: the fact that $R$ is grade-symmetric along with the Auslander condition tells us that $j_{R}(\mathrm{Ext}^{i}_{R}(R/P,R))\geq i\gneq j_{R}(R/P)$. As $-j_{R}$ is a dimension function, this implies that $\mathrm{Ext}^{i}_{R}(R/P,R)$ is torsion as an $R/P$-module. If $i\lneq j_{R}(R/P)$, then it is immediate from the definition that $\mathrm{Ext}^{i}_{R}(R/P,R)=0$
\end{proof}
\end{Lemma1}

If we instead choose to work with an Auslander dualizing complex $R$ over a pair of $k$-algebras $(A,B)$, similar results hold. There is an obvious analogue of \cref{extpure} (\cite[Theorem 2.14]{JY2}), and we can also recover an analogue of \cref{onlyone} if we replace $R$ with an Auslander, Cdim-symmetric, weakly bifinite dualizing complex.

\vspace{1.5mm}

We will need the following result concerning the purity of prime factor rings. It must be well-known, but we have been unable to find a suitable reference.

\begin{Lemma1} \label{primepure}
Let $R$ be an Auslander-Gorenstein ring. Let $P$ be a prime ideal. Then, $R/P$ is pure as a left and a right $R$-module.
\begin{proof}
As the hypothesis is symmetric it is enough to work over right modules. Let $K$ be a right ideal of $R/P$. It is immediate that $j_{R}(K)\geq j_{R}(R/P)$. However, as $R/P$ is right noetherian, $K$ contains a uniform right ideal $U$. If $n$ is the Goldie rank of $R/P$, then $R/P$ contains a copy of $\bigoplus U^{n}$ as an essential submodule (\cite[Proposition 7.23]{GW}), and thus there exists an embedding of $R/P$ into  $\bigoplus U^{n}$ - this is a consequence of Goldie's regular element lemma  (\cite[Proposition 6.13]{GW}). Since $U\leq K$, this implies that there is an embedding of $R/P$ into $\oplus K^{n}$. We conclude that $j_{R}(R/P)\geq j_{R}(K)$, and we have equality. Therefore every submodule of $R/P$ has the same grade as $R/P$, and $R/P$ is pure.
\end{proof}
\end{Lemma1}

\begin{Proposition1} \label{locduality}
Let $R$ be an Auslander-Gorenstein, grade-symmetric ring, and let $C$ be a localizable clique of prime ideals, with corresponding Ore set $X$. Let $P\in C$. Then, $R_{X}$ is an AS-Gorenstein ring of injective dimension $j_{R}(R/P)$, and there is a bijection $D:C\rightarrow C$ such that given $P\in C$ and any right module $M$ of finite flat dimension, there is an isomorphism: $$\mathrm{Tor}_{i}^{R_{X}}(M_{X}, {^{D(P)}S)}\cong \mathrm{Ext}^{j_{R}(R/P)-i}_{R_{X}}(S^{P},M_{X}).$$

\begin{proof}
We first show that $R_{X}$ is AS-Gorenstein.

 By \cref{primepure}, $R/P$ is pure as a right module for every prime ideal $P$. Non-zero localizations of pure modules are pure of the same grade. To see this, note that by \cite[Theorem 3.4.2.6]{z-f}, a module $M$ is pure if and only if $\mathrm{Ext}^{i}_{R}(\mathrm{Ext}^{i}_{R}(M,R),R)=0$ for all $i\neq j_{R}(M)$. Using \cref{flatten} twice, we see that $$\mathrm{Ext}^{i}_{R}(\mathrm{Ext}^{i}_{R}(M,R),R)\otimes_{R} R_{X}\cong \mathrm{Ext}^{i}_{R_{X}}(\mathrm{Ext}^{i}_{R_{X}}(M_{X},R_{X}),R_{X}).$$ Since $R/P$ is pure, by \cref{primepure}, we have that $R_{X}/P_{X}$ is pure, and $j_{R_{X}}(R_{X}/P_{X})=j_{R}(R/P)$ for every prime ideal in $C$. 

As $R$ is grade-symmetric, this implies that $j_{R_{X}}(R_{X}/P_{X})=j_{R_{X}}(R_{X}/Q_{X})$ for all $P$ and $Q$ in $C$. To see this, note that if $P\rightsquigarrow Q$, then there exists a bimodule $B$ which is finitely generated and torsionfree both as a left $R/P$-module and as a right $R/Q$-module. Thus $j_{R}(R/P)=j_{R}(_{R}B)$, and $j_{R}(R/Q)=j_{R}(B_{R})$. Since $j_{R}$ is a grade-symmetric dimension function, it follows that $j_{R}(R/P)=j_{R}(R/Q)$.

 However in $R_{X}$, $Q_{X}$, as $Q$ varies across $C$, are the only left and right primitive ideals of $R_{X}$, and they are all co-artinian.  Thus, every simple right or left module has grade equal to $j_{R}(R/P)$. The reader will easily check that localizations of Auslander-Gorenstein rings at Ore sets are Auslander-Gorenstein. For such a ring, there is always a simple module whose grade is equal to the injective dimension of the ring. Thus the injective dimension of $R_{X}$ is $j_{R}(R/P)$, and every simple module of $R_{X}$ is holonomic.

For an Auslander-Gorenstein ring $A$ of dimension $n$, $\mathrm{Ext}_{A}^{n}($---$,A)$ defines a duality of categories between holonomic left and right modules which preserves simple modules. This, along with the previous paragraph, implies that $R_{X}$ is AS-Gorenstein. 

We now proceed to establish the duality stated in the theorem. Given a prime ideal $P$ in $C$, we have that there is a unique prime ideal $D(P)\in C$ such that $D(P)_{X}$ annihilates $\mathrm{Ext}_{R_{X}}^{j_{R}(R/P)}(R_{X}/P_{X},R_{X})$ on the left; this establishes a bijection $D:C\rightarrow C$. In particular, the way we construct $D$ implies that $\mathrm{Ext}^{j_{R}(R/P)}_{R_{X}}(S^{P},R_{X})\cong  \mathrm{^{D(P)}S}.$

Let $M$ now be a right $R$-module of finite flat dimension. Then, by \cref{iss1}, there is a convergent spectral sequence $$E^{2}_{p,q}= \mathrm{Tor}_{p}^{R_{X}}(M_{X}, \mathrm{Ext}_{R_{X}}^{q}(S^{P},R_{X}))\Rightarrow \mathrm{Ext}_{R_{X}}^{q-p}(S^{P},R_{X}).$$ Since every simple module is holonomic, this collapses to give us the proposed isomorphism.
\end{proof}
\end{Proposition1}

\subsection*{Depth}

We define the \textit{depth} of a module $M$ over a ring $R$ to be $$\depth{M}=\mathrm{inf}\{i|\exists S\ simple\  s.t.\  \mathrm{Ext}^{i}_{R}(S,M)\neq 0\}.$$  

In \cite{QZ1}, Wu and Zhang prove the following result:

\begin{Theorems1} \cite[Theorem 0.1]{QZ1} \label{niceclique2}
Let $R$ be an AS-Gorenstein algebra over a field $k$ of injective dimension $d$ such that $R$ is either FBN or semilocal. Let $M$ be a finitely generated non-zero right $R$-module.
\begin{itemize}
 \item [i.)](Auslander-Buchsbaum formula) If the projective dimension of $M$, $\pd{M}$, is finite, then $\pd{M}+\depth{M}=d$.
 \item [ii.)](Bass's Theorem) If the injective dimension of $M$, $\id{M}$, is finite, then $\id{M}=d$.
 \item [iii.)](No-holes theorem) For every $i$ such that $\depth{M}\leq i\leq \id{M}$, there exists a simple $R$-module $S$ such that $\mathrm{Ext}^{i}_{R}(S,M)\neq 0$.
\end{itemize}
\end{Theorems1}

We will prove that the Auslander-Buchsbaum formula, Bass's Theorem, and the No-holes theorem hold for a class of rings including that include the localizations of grade-symmetric Auslander-Gorenstein algebras over uncountable fields at localizable cliques. The following result is essentially an application of the spectral sequence in \cref{iss1} to \cite[Theorem 12]{W}. We refer the reader to \cite[Chapter 16]{GW} for details about the generic regularity condition. The reader should note that collections of prime ideals satisfying the generic regularity condition are not rare; if $X$ is a set of prime ideals in a noetherian ring $R$ such that $\{\rank{P}|P\in X\}$ is finite, then $X$ satisfies the generic regularity condition. In particular, every finite set of prime ideals in a noetherian ring satisfies the generic regularity condition, as does every set of completely prime ideals in a noetherian ring. 

\begin{Proposition1} \label{niceclique}
Let $R$ be an AS-Gorenstein ring of dimension $d$ which is either FBN or satisfies the following conditions:
\begin{itemize}
\item[a.)] Every right or left primitive ideal is co-artinian.
\item[b.)] The set of co-artinian ideals is finite, or $R$ contains a set of central units $F$ such that the difference of two distinct elements of $F$ is a unit, and the cardinality of $F$ is larger than the cardinality of the set of co-artinian ideals.
\item[c.)] The set of co-artinian prime ideals satisfies the generic regularity condition.
\end{itemize}
 Let $M\neq 0$ be a finitely generated module over $R$, such that the flat dimension of $M$ is finite. Then,
\begin{itemize}
 \item [i.)] $\pd{M}+\depth{M}=d$.
 \item [ii.)] $\id{M}=d$.
 \item [iii.)] For every $i$ such that $\depth{M}\leq i\leq \id{M}$, there exists a simple $R$-module $S$ such that $\mathrm{Ext}^{i}_{R}(S,M)\neq 0$.
\end{itemize}
\begin{proof}

Over a ring which satisfies these hypotheses, the right primitive ideals coincide with the maximal ideals and the left primitive ideals, and happen to be co-artinian. This follows from \cite[Proposition 9.4]{GW} for FBN rings, and follows from the hypotheses in the other case. Thus every right or left simple module $S$ is a summand of $R/P$ for some maximal ideal $P$, viewed as either a right or a left module. For the purpose of this proof, let $^{P}S$ (resp. $S^{P}$) denote the unique left (resp. right) simple module with annihilator $P$.

We can now use the argument described at the beginning of \cite[Theorem 12]{W}: there is a spectral sequence which  collapses to give us the following isomorphism: $$\mathrm{Hom}_{R/P}(\mathrm{Tor}_{n}^{R}(L,R/P),S^{P})\cong \mathrm{Ext}^{n}_{R}(L,S^{P}).$$ It follows that $\mathrm{Ext}^{p}_{R}(L,S^{P})\neq 0$ if and only if $\mathrm{Tor}_{R}^{p}(L,^{P}S)\neq 0$ for any right module $L$.

If $M$ is a finitely generated module such that $\mathrm{Tor}_{1}^{R}(M,S)=0$ for all simple modules $S$, then $M$ is projective. This is a consequence of \cite[Lemma 1.1]{QZ1} for FBN rings, and from \cite[Theorem 12]{W} for rings satisfying the hypotheses for the other case of the proposition. While \cite[Lemma 1.1]{QZ1} is stated for algebras defined over a field, the proof continues to hold without any modification for general rings. By induction, for every $i$ between $0$ and the projective dimension of $M$, there exists a simple module $S$ such that $\mathrm{Tor}_{i}^{R}(M,S)\neq0$. 

We now use the fact that our rings are AS-Gorenstein. By an argument very similar to that used in \cref{locduality}, there is a bijection $D$ with inverse $D'$ between the right and left simple modules such that $\mathrm{Ext}_{R}^{i}(S,M)\cong \mathrm{Tor}_{d-i}^{R}(M,D(S))$. If $S$ is a simple factor of $M$ with annihilator $P$, then $M\otimes_{R} R/P \neq 0$ and therefore $\mathrm{Ext}_{R}^{d}(D'(S),M)\neq 0$. The result now follows.
\end{proof}
\end{Proposition1}

\begin{Rem}
The reader may wonder why we have included FBN rings in the statement of \cref{niceclique} when an ostensibly stronger result is proved in \cref{niceclique2}. The techniques involved in the proof of \cref{niceclique2} involve the use of dualizing complexes, and the theory of dualizing complexes can be quite delicate when working with rings which are not algebras over a field. Since we want to consider FBN rings which are not algebras over a field, we include them in the statement of \cref{niceclique}.
\end{Rem}

\section{Injective Modules, the Second Layer Condition, and FBN Rings}

In this section, we make our first attempt at localizing $\mathrm{Ext}$.

\begin{Lemma1} \label{localniceev}
Let $R$ be a right noetherian ring, and let $X$ be a right denominator set. Let $M$ be a right $R$-module with an injective resolution $\mathrm{\mathbf{I}} := 0\rightarrow I^{0}\rightarrow I^{1}\rightarrow...$ such that every indecomposable direct summand of every $I^{i}$ is either $X$-torsion or $X$-torsionfree, for all $i$. Let $P$ be a prime ideal in $R$. Then, for all $i$, there is an isomorphism $$\mathrm{Ext}^{i}_{R}(R/P,M)\otimes_{R} R_{X}\cong \mathrm{Ext}^{i}_{R_{X}}(R_{X}/P_{X},M_{X}).$$
If, in addition, $S$ is a ring and $M$ is an $(S,R)$-bimodule, the above isomorphism is one of bimodules. 
\begin{proof}
Since $R_{X}$ is flat as a left $R$-module, it follows that $\mathrm{Ext}^{i}_{R}(R/P,M)\otimes_{R} R_{X}$ is the $i^{th}$ homology of the complex 
\begin{equation} \label{homo1}
\mathrm{Hom}_{R}(R/P,\mathrm{\mathbf{I}})\otimes_{R} R_{X}\cong \mathrm{ann}_{\mathrm{\mathbf{I}}}(P)\otimes_{R} R_{X}.
\end{equation}

 If $E$ is an $X$-torsion module, then ${E}_{X}$ and $\mathrm{ann}_{E}(P)_{X}$ are both $0$. However, if $E$ is $X$-torsionfree, it follows that $\mathrm{ann}_{E}(P)_{X}\cong \mathrm{ann}_{E_{X}}(P_{X})$. 

Moreover, if ${E}$ is an injective $R$-module which is either $X$-torsion or $X$-torsionfree, then ${E_{X}}$ is an injective $R_{X}$-module. More precisely, if ${E}$ is $X$-torsion, then ${E_{X}}=0$, while if it is $X$-torsionfree, then by \cite[Corollary 10.14]{GW}, ${E}$ carries an $R_{X}$-module structure and is injective as an $R_{X}$-module; in particular, ${E_{X}}$ is injective. 

From this, it follows that $\mathrm{\mathbf{I}}_{X}$ is an injective resolution of $M_{X}$ over $R_{X}$, and thus that $\mathrm{Ext}^{i}_{R_{X}}(R_{X}/P_{X},M_{X})$ is the homology of the complex $$\mathrm{Hom}_{R_{X}}(R_{X}/P_{X},\mathrm{\mathbf{I}}_{X})\cong \mathrm{ann}_{\mathrm{\mathbf{I}}_{X}}(P_{X}).$$ However, by what we have shown in the first paragraph of this argument, this complex is isomorphic to the complex displayed in \ref{homo1}.  In particular, they have isomorphic homologies, and we conclude that for all $i$, there is an isomorphism as in the statement of the lemma.

If $M$ happens to be an $(S,R)$-bimodule, then there is an action of $S$ on the injective resolution $\mathrm{\mathbf{I}}$ in the homotopy category, which comes from lifting the action on $M$. This gives the homology of $\mathrm{\mathbf{I}}$, and thus the various derived functors we are constructing a left $S$-module structure. It is now immediate that the isomorphism we have constructed above preserves this structure. 
\end{proof}
\end{Lemma1}

In \cite[Theorem 1.1]{GJ} (a result attributed to Brown), Goodearl and Jordan prove that over a noetherian right FBN ring, given a denominator set $X$, indecomposable injective modules are either $X$-torsion or $X$-torsionfree. In \cite[Theorem 1]{GJ2}, Goodearl and Jordan characterise the denominator sets with this property: if $X$ is a right denominator set in a right noetherian ring $R$, every indecomposable injective module is either $X$-torsion or $X$-torsionfree if and only if every essential extension of a $X$-torsion module is $X$-torsion if and only if localization at $X$ preserves essential extensions. It is worth noting that the Ore set associated to a localizable clique need not preserve essentiality under localization (see \cite[Example 6]{GJ2}).

\begin{Coroll1}[of \cref{localniceev}]
Let $R$ be a right noetherian ring, and let $X$ be a right denominator set such that localization at $X$ preserves essential extensions. Let $P$ be a prime ideal in $R$. Let $M$ be any $R$-module. Then, for all $i$, there is an isomorphism $$\mathrm{Ext}^{i}_{R}(R/P,M)\otimes_{R} R_{X}\cong \mathrm{Ext}^{i}_{R_{X}}(R_{X}/P_{X},M_{X}).$$
If, in addition, $S$ is a ring and $M$ is an $(S,R)$-bimodule, then the above isomorphism is one of bimodules. 
\end{Coroll1}

We now prove a proposition in this section deals with the behaviour of tame injective modules associated to suitably nice prime ideals when localizing. This is a minor generalization of the argument used to prove \cite[Theorem 1.1]{GJ}. Indeed, the argument that we use is essentially the argument in \cite[Theorem 1.1]{GJ}, but in a more modern language which allows us to see the generalization.

\begin{Proposition1} \label{localizeinj}
Let $R$ be a noetherian ring, and let $P$ be a prime ideal such that every prime in the right link closure of $P$ satisfies the right second layer condition. Then, given a denominator set $X$, we have that $({E_{P}})_{X}$ is either $X$-torsion or $X$-torsionfree. If $P\cap X = \emptyset$, then $({E_{P}})_{X}\cong E_{P_{X}}$, otherwise $({E_{P}})_{X}=0$.
\begin{proof}
We will show that $E_{P}$ is either $X$-torsion or $X$-torsionfree, depending on whether $P\cap X \neq \emptyset$ or not. In the first case, $E_{P}$ vanishes when localized; in the second, it itself carries an $R_{X}$-module structure and is injective as an $R_{X}$-module, by \cite[Corollary 10.14]{GW}.

Suppose $P\cap X = \emptyset$. Then, \cite[Lemma 10.19]{GW} informs us that $X\subset C_{R}(P)$, the set of regular elements modulo $P$. We claim that $E_{P}$ is $X$-torsionfree. If not, then the $X$-torsion submodule $t_{X}(E_{P}) \neq 0$, and thus contains a submodule isomorphic to a uniform right ideal of $R/P$, say $U$. This is a contradiction, as $X\subset C_{R}(P)$. Therefore $P\cap X = \emptyset$ implies that $E_{P}$ is $X$-torsionfree. It is then immediate that as a $R_{X}$-module, it is isomorphic to $E_{P_{X}}$.

Now let us assume that $P\cap X \neq \emptyset$. We will show that every finitely generated submodule of $E_{P}$ is $X$-torsion. Let $M$ be such a submodule. As every prime in the right link closure of $P$ satisfies the right second layer condition, there is an affiliated series for $M$, $0 \leq M_{1} \leq M_{2} \leq ... \leq M_{n}= M$ with corresponding affiliated primes $P_{i}$, such that $M_{i+1}/M_{i}$ is isomorphic to a torsionfree module over $R/P_{i}$. If $M$ is not $X$-torsion, then  some $M_{i+1}/M_{i}$ is not $X$-torsion, which then implies that $X$ intersects $P_{i}$ trivially. However, as $X$ is also a left denominator set and $P_{i}$ lies in the right link closure of $P$,  \cite[Lemma 14.7]{GW} implies that $X$ intersects $P$ trivially, a contradiction. Thus $E_{P}$ is $X$-torsion.
\end{proof}
\end{Proposition1}

\begin{Coroll1} \label{localizeext}
Let $R$ be a noetherian ring which satisfies the right second layer condition. Let $X$ be an Ore set in $R$. Let $M$ be an $R$-module such that the minimal injective resolution of $M$ contains only tame injective modules. Then, given a prime ideal $P$, there is an isomorphism of $R_{X}$-modules for all $i$: $$ \mathrm{Ext}^{i}_{R}(R/P,M)\otimes_{R} R_{X}\cong \mathrm{Ext}^{i}_{R_{X}}(R_{X}/P_{X},M_{X}).$$ If, in addition, $M$ is an $(S,R)$-bimodule, then the above isomorphism is one of bimodules. 
\begin{proof}
This follows from \cref{localniceev} and \cref{localizeinj}.
\end{proof}
\end{Coroll1}

\begin{Rem}
Modules which satisfy the hypotheses of \cref{localizeext} occur with greater frequency than one might originally suspect. For example, let $R$ be a ring which satisfies the second layer condition. If $M$ is an $R$-module such that for every associated prime $P$ of $M$, $R/Q$ is FBN for every prime in the clique of $P$, it is an easy exercise to show that $M$ satisfies the hypotheses of \cref{localizeext}. For instance, co-artinian primes ideals are good examples of such prime ideals. 
\end{Rem}

\section{The bimodule $\mathrm{Ext}^{j_{R}(R/P)}_{R}(R/P,R)$}

In this section, we study the structure of $\mathrm{Ext}^{j_{R}(R/P)}_{R}(R/P,R)$ as an $(R,R)$-bimodule. Our basic result has the following form: given a prime ideal $P$, there exists a prime ideal $Q$ such that  $\mathrm{Ext}^{j_{R}(R/P)}_{R}(R/P,R)$ is torsionfree as a left module over $R/Q$ and torsionfree as a right module over $R/P$. We will give two proofs of our result.
In \textsection 3.1, we give an argument that depends heavily on localization; this proof is somewhat elementary, but only works for a reasonably restrictive class of rings.
In \textsection 3.2, we elaborate on an observation of James Zhang that the theory of dualizing complexes gives us a similar result for a much larger collection of rings. The reader should note the power and elegance that comes with using the theory of dualizing complexes. 

\subsection*{The First Proof: Exploiting a Local Duality}

Let $E$ be an injective left module. For all finitely generated right modules $M$, there is an isomorphism: 

\begin{equation} \label{niceice}
\mathrm{Hom}_{R}(\mathrm{Ext}^{i}_{R}(M,R),E)\cong \mathrm{Tor}_{i}^{R}(M,E).
\end{equation}

This formula follows immediately from the collapsing of the spectral sequence in \cref{iss2}.

If $M$ happens to be an $(S,R)$-bimodule, then the groups in the above formula are left $S$-modules, and the above isomorphism preserves this structure. 

The following easy observation is key:

\begin{Lemma1} \label{localhom}
Let $R$ and $S$ be rings, and suppose $M$ is an $(R,S)$-bimodule which has finite uniform rank as an $R$-module and is torsionfree as an $S$-module. Let $N$ be a left $R$-module, and suppose there exists a homomorphism $f:M\rightarrow N$ such that $\ker(f)$ is not essential as an $R$-module. Then, $\mathrm{Hom}_{R}(M,N)$ is not torsion as a left $S$-module.
\begin{proof}
If there exists a $c\in C_{S}(0)$ such that $cf=0$, then $Mc\subset \ker(f)$. However, as $M$ is torsionfree, $Mc\cong M$, and the finite uniform rank of $M$ forces $Mc$ to be essential in $M$. Thus $\ker(f)$ is essential in $M$, a contradiction. Therefore $f$ is not annihilated by a regular element, and thus is not torsion.
\end{proof}
\end{Lemma1}

\begin{Proposition1} \label{ass1}
Let $R$ be a noetherian ring, and let $P$ be a prime ideal such that $\mathrm{Ext}^{i}_{R}(R/P,R)$ is torsionfree as a right $R/P$-module. Then, $Q$ is an associated prime of $\mathrm{Ext}^{i}_{R}(R/P,R)$ as a left module only if there exists some injective module $E$ with unique associated prime $Q$ such that $\mathrm{Tor}^{R}_{i}(\mathcal{Q}(R/P),E)\neq 0$. If $\mathrm{Ext}^{i}_{R}(R/P,R)$ is also finitely generated as a right $R/P$-module, then $Q$ is an associated prime only if  $\mathrm{Tor}^{R}_{i}(\mathcal{Q}(R/P),E_{Q})\neq 0$.
\begin{proof}

If $Q$ is an associated prime of  $\mathrm{Ext}^{i}_{R}(R/P,R)$ as a left module, then there exists a submodule which is annihilated by $Q$ and is fully faithful as an $R/Q$-module. Thus, $\mathrm{Ext}^{i}_{R}(R/P,R)$, as a left module, has a uniform submodule which is fully faithful as an $R/Q$-module. Let the injective hull of this submodule be $E$. By the injectivity of $E$, there is a homomorphism from $\mathrm{Ext}^{i}_{R}(R/P,R)$ to $E$ such that the kernel of this map is not essential. Thus, by \cref{localhom},  $\mathrm{Hom}_{R}(\mathrm{Ext}^{i}_{R}(R/P,R),E)$ is not torsion as an $R/P$-module. 

It follows from \ref{niceice} for $M=R/P$ that there is an isomorphism $$\mathrm{Hom}_{R}(\mathrm{Ext}^{i}_{R}(R/P,R),E)\cong \mathrm{Tor}^{R}_{i}(R/P,E).$$ Thus $\mathrm{Tor}^{R}_{i}(R/P,E)$ is not torsion as an $R/P$-module. Tensoring with $\mathcal{Q}(R/P)$, we see that  $\mathcal{Q}(R/P)\otimes_{R} \mathrm{Tor}^{R}_{i}(R/P,E)\neq 0$. However, $$\mathcal{Q}(R/P)\otimes_{R} \mathrm{Tor}^{R}_{i}(R/P,E)\cong \mathrm{Tor}^{R}_{i}(\mathcal{Q}(R/P),E).$$

If, in addition, we have that $\mathrm{Ext}^{i}_{R}(R/P,R)$ is finitely generated as a $R/P$-module, then given an associated prime $Q$, we have that the affiliated submodule corresponding to $Q$ is torsionfree as a $R/P$-module by \cite[Lemma 8.1]{GW}. Thus, we can take $E$ to be $E_{Q}$, the unique tame indecomposable injective with assassinator $Q$. The argument then follows as before. 
\end{proof}
\end{Proposition1}

\begin{Theorems1} \label{bimodstruc}
Let $R$ be an Auslander-Gorenstein, grade-symmetric, weakly bifinite ring which satisfies the second layer condition. Let $P$ be a prime ideal such that the clique of $P$, $C$, is localizable, with corresponding Ore set $X$. Then, there is a bijection $D:C\rightarrow C$ such that $\mathrm{Ext}^{j_{R}(R/P)}_{R}(R/P,R)$ is annihilated as a left module by $D(P)$, and is torsionfree both as a right $R/P$-module and a left $R/D(P)$ module.
\begin{proof}
The map $D$ is the map described in \cref{locduality}.

We begin by showing that $\mathrm{Ext}^{j_{R}(R/P)}_{R}(R/P,R)$ has a unique associated prime as a left $R$-module. Suppose $Q$ is an associated prime of $\mathrm{Ext}^{j_{R}(R/P)}_{R}(R/P,R)$. By \cref{ass1}, we know that $\mathrm{Tor}_{j_{R}(R/P)}^{R}(\mathcal{Q}(R/P),E_{Q})\neq 0$. 

There is an isomorphism $$\mathrm{Tor}_{j_{R}(R/P)}^{R}(\mathcal{Q}(R/P),E_{Q})\cong \mathrm{Tor}_{j_{R}(R/P)}^{R_{X}}(R_{X}/P_{X}, R_{X}\otimes_{R} E_{Q}).$$ This follows from a direct computation, after noting that $R_{X}$ is flat both as a left and as a  right $R$-module. Since  $\mathrm{Tor}_{j_{R}(R/P)}^{R}(\mathcal{Q}(R/P),E_{Q})\neq 0$, we conclude using \cref{localizeinj} that $Q\cap X = \emptyset$, and $R_{X}\otimes_{R} E_{Q}$ is isomorphic to the tame indecomposable injective $E_{Q_{X}}$. 

We continue to use the notation from \cref{locduality}. By \cref{locduality},  $$\mathrm{Tor}_{j_{R}(R/P)}^{R_{X}}(R_{X}/P_{X}, E_{Q_{X}})\cong \mathrm{Hom}_{R_{X}}(^{D(P)}S,E_{Q_{X}})^{l(\mathcal{Q}(R/P))}.$$ Thus there is a non-zero homomorphism $^{D(P)}S\rightarrow E_{Q_{X}}$, which by the simplicity of $^{D(P)}S$ must be injective. This fact, along with the maximality of $D(P)_{X}$ in $R_{X}$, implies that $Q_{X}=D(P)_{X}$, and so $Q=D(P)$. We have thus proved that  $\mathrm{Ext}^{j_{R}(R/P)}_{R}(R/P,R)$  has a unique associated prime, $D(P)$, which by construction lies in the clique of $P$. We now show that  $\mathrm{Ext}^{j_{R}(R/P)}_{R}(R/P,R)$ is annihilated by $D(P)$ on the left, and is torsionfree as an $R/D(P)$ module. 

Tensoring by $R_{X}$ on the left, \cref{flatten} implies that $$R_{X}\otimes_{R} \mathrm{Ext}^{j_{R}(R/P)}_{R}(R/P,R)\cong \mathrm{Ext}^{j_{R}(R/P)}_{R_{X}}(R_{X}/P_{X},R_{X}).$$ By the proof of \cref{locduality}, the latter bimodule is an $(R_{X}/D(P)_{X},R_{X}/P_{X})$-bimodule, and thus satisfies the conclusions of our theorem i.e. it is annihilated by $D(P)$ on the left, is torsionfree as a left module over $R/D(P)$, and is a torsionfree right $R/P$-module. If we let $W$ denote the left $X$-torsion bimodule of $\mathrm{Ext}^{j_{R}(R/P)}_{R}(R/P,R)$, then there is a natural embedding $\mathrm{Ext}^{j_{R}(R/P)}_{R}(R/P,R)/W\rightarrow \mathrm{Ext}^{j_{R}(R/P)}_{R_{X}}(R_{X}/P_{X},R_{X})$. If we can show that $W=0$, then we are done.

If $W\neq 0$, then it has an associated prime, which must be $D(P)$. However, this implies that $W$ contains a fully faithful torsion $R/D(P)$ module, as $W$ is $X$-torsion, and $D(P)$ belongs to the clique $C$. Let this module be denoted by $V$. Then $j_{R}(V)\gneq j_{R}(R/D(P))=j_{R}(R/P)$, as $-j_{R}$ is a dimension function. As $\mathrm{Ext}^{j_{R}(R/P)}_{R}(R/P,R)$ is pure as a left $R$-module, by \cref{extpure}, this is a contradiction.  Thus  $W=0$, and we have proved our result.
\end{proof}
\end{Theorems1}

\begin{Rem}
We would like to prove the above result by only assuming that $R$ satisfies the second layer condition locally, i.e. by only assuming that the clique of $P$ satisfies the second layer condition. However, our arguments do not work in this generality; we need the second layer condition to hold globally in order to understand the structure of the injective module $E$ after localizing at $X$. In the next section, we show that by using the theory of dualizing complexes, the result continues to hold in certain situations if we only assume a local hypothesis. 
\end{Rem}

\begin{Coroll1} \label{bimod2}

Let the hypothesis of \cref{bimodstruc} be satisfied. Then, there are isomorphisms of $(R,R)$-bimodules: 

\begin{align*}
\mathrm{Ext}^{j_{R}(R/P)}_{R}(R/P,R)\otimes_{R} R_{X} & \cong R_{X}\otimes_{R} \mathrm{Ext}^{j_{R}(R/P)}_{R}(R/P,R) \\
& \cong   \mathrm{Ext}^{j_{R}(R/P)}_{R_{X}}(R_{X}/P_{X},R_{X}).\\
\end{align*}

\begin{proof}

The fact that $$R_{X}\otimes_{R} \mathrm{Ext}^{j_{R}(R/P)}_{R}(R/P,R)\cong \mathrm{Ext}^{j_{R}(R/P)}_{R_{X}}(R_{X}/P_{X},R_{X})$$ is a standard application of \cref{flatten}. It follows that $R_{X}\otimes_{R} \mathrm{Ext}^{j_{R}(R/P)}_{R}(R/P,R)$ carries the structure of a right $R_{X}/P_{X}$ module; in particular, it also carries the structure of an $R_{X}$-module. Therefore there is an isomorphism $$R_{X}\otimes_{R} \mathrm{Ext}^{j_{R}(R/P)}_{R}(R/P,R)\cong R_{X}\otimes_{R} \mathrm{Ext}^{j_{R}(R/P)}_{R}(R/P,R)\otimes_{R} R_{X}.$$ 

The remainder of the argument is similar to the one used in \cite[Proposition 3.3]{JY}. By \cref{bimodstruc}, $\mathrm{Ext}^{j_{R}(R/P)}_{R}(R/P,R)$ is torsionfree as a right $R/P$-module. Tensoring with $R_{X}$ on the right, we immediately see that  $\mathrm{Ext}^{j_{R}(R/P)}_{R}(R/P,R)$ embeds as an essential subbimodule of $\mathrm{Ext}^{j_{R}(R/P)}_{R}(R/P,R)\otimes_{R} R_{X}.$ Thus $\mathrm{Ext}^{j_{R}(R/P)}_{R}(R/P,R)\otimes_{R} R_{X}$ is torsionfree as a left $R/D(P)$ module (essential extensions of torsionfree modules are torsionfree), and the fact that it is finitely generated as a right $\mathcal{Q}(R/P)$ module implies that it  carries a left $\mathcal{Q}(R/D(P))$ module structure - this last statement is true because left multiplication by any element $x\in C_{R}(D(P))$ acts as an injective right $\mathcal{Q}(R/P)$ homomorphism; the fact that the module is finite length as a $\mathcal{Q}(R/P)$ module tells us that this left multiplication map must be surjective. Therefore there is an isomorphism $$\mathrm{Ext}^{j_{R}(R/P)}_{R}(R/P,R)\otimes_{R} R_{X} \cong R_{X}\otimes_{R} \mathrm{Ext}^{j_{R}(R/P)}_{R}(R/P,R)\otimes_{R} R_{X}.$$ 

Thus $$\mathrm{Ext}^{j_{R}(R/P)}_{R}(R/P,R)\otimes_{R} R_{X}\cong R_{X}\otimes_{R} \mathrm{Ext}^{j_{R}(R/P)}_{R}(R/P,R)$$ as $(R,R)$-bimodules.
\end{proof}
\end{Coroll1}

\subsection*{The Second Proof: The Approach Through Dualizing Complexes}

As mentioned earlier, the author is indebted to James Zhang for pointing out the approach described in this section.  Most of this section is a description of results in the literature that we then piece together. 

Most of the literature about dualizing complexes has been about dualizing complexes defined over fields. However, we will work in a more general setting: in this subsection, we fix a commutative base ring $k$, and demand that all of our rings are $k$-algebras that are flat over $k$. The definition of a dualizing complex as it is defined for algebras over a field seems to be the correct definition of a dualizing complex in this situation as well. 

As a convention, whenever we consider a dualizing complex over the pair $(A,B)$, we assume that $A$ is left noetherian and $B$ is right noetherian, and that derived functors and other constructions over $A$ are defined using left modules, while those defined over $B$ use right modules. 

We begin by describing the notion of an induced dualizing complex: 

\begin{Definitions1} \cite[Definition 4.1]{JY}
Let $R$ be a dualizing complex over $(A,B)$. Let $A\rightarrow \bar{A}$ and $B\rightarrow \bar{B}$ be algebra homomorphisms such that $\bar{A}$ is a finite left $A$-module and $\bar{B}$ is a finite right $B$ module. Then, a dualizing complex $\bar{R}$ over $(\bar{A},\bar{B})$ is \textit{induced} by $R$ if there are morphisms, called the left and right trace morphisms, $$Tr_{l}, Tr_{r}:\bar{R}\rightarrow R$$ in $D^{b}({A^{op}\otimes_{k} B})$ such that the following induced maps are isomorphisms in  $D^{b}({\bar{A}^{op}\otimes_{k} B})$ and  $D^{b}({A^{op}\otimes_{k}\bar{B}})$ respectively: $$\bar{R}\rightarrow \mathrm{RHom}_{A}(\bar{A},R),$$ $$\bar{R}\rightarrow \mathrm{RHom}_{B}(\bar{B},R).$$
\end{Definitions1}

\cite[Lemma 4.3]{JY} implies that any dualizing complex induced by an Auslander dualizing complex is Auslander.

The following definition is a minor extension of \cite[Definition 4.4]{JY}: 

\begin{Definitions1} 
Let $\mathcal{F}_{A}$ be a class of ideals of $A$ and $\mathcal{F}_{B}$ a class of ideals of $B$. We say that a dualizing complex $R$ over $(A,B)$ has the \textit{trace property between $\mathcal{F}_{A}$ and $\mathcal{F}_{B}$} if there exists an inclusion preserving bijection $\phi$ between $\mathcal{F}_{B}$ and $\mathcal{F}_{A}$ such that given an ideal $I\in \mathcal{F}_{B}$, there is a dualizing complex $\bar{R}$ over $(A/\phi(I), B/I)$ that is induced by $R$.
\end{Definitions1}

In particular, following \cite[Definition 4.4]{JY} we will say that a dualizing complex $R$ over $(A,B)$ has the \textit{trace property for ideals} when it has the trace property between $\mathcal{F}_{A}$ and $\mathcal{F}_{B}$ where $\mathcal{F}_{A}$ is the collection of all the ideals of $A$, and $\mathcal{F}_{B}$ is the collection of all the ideals of $B$. We say that a dualizing complex $R$ over $(A,B)$ has the \textit{trace property for prime ideals} when it has the trace property between $\mathcal{F}_{A}$ and $\mathcal{F}_{B}$ where $\mathcal{F}_{A}$ is the spectrum of $A$ and $\mathcal{F}_{B}$ is the spectrum of $B$. 

\cite[Lemma 4.6]{JY} tells us that if an Auslander dualizing complex has the trace property between $\mathcal{F}_{A}$ and $\mathcal{F}_{B}$ with respect to a bijection $\phi:\mathcal{F}_{B} \rightarrow \mathcal{F}_{A}$, then $\phi$ induces a bijection between the prime ideals in $\mathcal{F}_{B}$ and the prime ideals in $\mathcal{F}_{A}$. In particular, if an Auslander dualizing complex has the trace property for ideals, it has the trace property for prime ideals.

We note here that \cite[Proposition 4.5]{JY} gives us a list of dualizing complexes that are known to satisfy the trace property for prime ideals. 

The following result should be compared to \cref{bimodstruc} and \cref{bimod2}:

\begin{Proposition1} \label{dual1}
Let $R$ be an Auslander dualizing complex over a pair of algebras $(A,B)$ such that $R$ satisfies the trace property between classes of prime ideals $\mathcal{F}_{A}$ and $\mathcal{F}_{B}$ with respect to a bijection $\phi: \mathcal{F}_{B}\rightarrow \mathcal{F}_{A}$. Let $P$ be a prime ideal in $\mathcal{F}_{B}$. Then, for all $i$, $\mathrm{Ext}^{i}_{B}(B/P,R)$ is an $(A/\phi(P),B/P)$-bimodule which is finitely generated on both sides, and $\mathrm{Ext}^{j_{R,B}(B/P)}_{B}(B/P,R)$ is torsionfree on both sides. There is an isomorphism $$\mathcal{Q}(A/\phi(P))\otimes_{A/\phi(P)} \mathrm{Ext}^{j_{R,B}(B/P)}_{B}(B/P,R)\cong \mathrm{Ext}^{j_{R,B}(B/P)}_{B}(B/P,R)\otimes_{B/P} \mathcal{Q}(B/P)$$

\begin{proof}
The argument essentially follows from the definitions and \cite[Proposition 3.3]{JY} and \cite[Lemma 4.3]{JY}. Note that the proofs of both of these results are still valid if $A$ and $B$ are flat algebras over a commutative ring $k$. We indulge ourselves by going through the details.

Suppose $R$ is a dualizing complex over $(A,B)$ which satisfies the trace property between $\mathcal{F}_{A}$ and $\mathcal{F}_{B}$. Then, by definition, given a prime ideal $P$ in $\mathcal{F}_{B}$, there exists a prime ideal $\phi(P)$ in $A$ and a dualizing complex $\bar{R}$ over $(A/\phi(P), B/P)$ induced by $R$. In particular, the cohomology of $\mathrm{RHom}_{B}(B/P,R)$ is annihilated by $\phi(P)$ on the left, and therefore the modules $\mathrm{Ext}^{i}_{B}(B/P,R)$ are $(A/\phi(P),B/P)$-bimodules for all $i$. As $\mathrm{RHom}_{B}(B/P,R)$ is a dualizing complex over $(A/\phi(P),B/P)$, by definition $\mathrm{Ext}^{i}_{B}(B/P,R)$ will be finitely generated on both sides.

Since a dualizing complex induced by an Auslander dualizing complex is Auslander, it follows that $\mathrm{RHom}_{B}(B/P,R)$ is an Auslander dualizing complex over $(A/\phi(P),B/P$), which are both prime noetherian rings. The degree of the lowest non-zero cohomology of $\mathrm{RHom}_{B}(B/P,R)$ is by definition $j_{R,B}(B/P)$ (which is equal to $j_{R,A}(A/\phi(P))$). \cite[Proposition 3.3.(1)]{JY}, and its proof, now complete the argument.
\end{proof}
\end{Proposition1}

As mentioned in the beginning of this section, our interest is in understanding the structure of $\mathrm{Ext}^{j_{A}(A/P)}_{A}(A/P,A)$ as a bimodule. Thus we will be mainly interested in applying \cref{dual1} to the situation where the dualizing complex $R$ is of the form $^{\lambda}A^{\rho}[n]$, where $\lambda$ and $\rho$ are some $k$-algebra automorphisms of $A$. Recall from \cref{elemhomoiden} that as a $(A,A)$-bimodule, $^{\lambda}A^{\rho}\cong {^{\rho^{-1} \lambda}A^{1}}$, via the map $\rho$. Thus we can restrict ourselves to twists on one side. This situation occurs with considerable frequency: the class of rings which have an Auslander, Cdim-symmetric dualizing complex of the form $^{\phi}A^{1}$ with the trace property for prime ideals include $k$-algebras $A$ with a noetherian connected filtration such that $gr(A)$ is AS-Gorenstein, and is either a.) PI, b.)  graded-FBN, or c.) has enough normal elements (see \cite[Lemma 2.1]{JY}, \cite[Proposition 4.5]{JY}, \cite[Corollary 6.9]{JY2}, and \cite[Proposition 6.18]{JY2}). It also includes complete, Auslander-Gorenstein, local $k$-algebras $A$, where $k$ is a field, with maximal ideal $J$ such that $A/J$ has finite dimension over $k$  (see \cite[Theorem 0.1]{CWZ}, \cite[Lemma 3.8]{CWZ}, and \cite[Proposition 4.5]{JY}). 

\begin{Coroll1} \label{dual2}
Let $A$ be a $k$-algebra such that there is an Auslander dualizing complex $R={^{\phi}A^{1}}[n]$ for some $\phi \in \mathrm{Aut(A)}$ which satisfies the trace property between classes of prime ideals $\mathcal{F}_{A}$ and $\mathcal{F'}_{A}$ with respect to a bijection $\psi: \mathcal{F}_{A}\rightarrow \mathcal{F'}_{A}$. Then, for all prime ideals $P$ in $\mathcal{F}_{A}$, $\mathrm{Ext}^{j_{A}(A/P)}_{A}(A/P,A)$ is a $(A/\phi(\psi(P)),A/P)$-bimodule which is torsionfree and finitely generated on both sides.
\begin{proof}
 Note that the structure of our dualizing complex implies that $A$ is Auslander-Gorenstein, and that as bimodules, $$\mathrm{Ext}^{j_{A}(A/P)}_{A}(A/P,{^{\phi}A^{1}})\cong \mathrm{Ext}^{j_{R}(A/P)}_{A}(A/P,R).$$ Now, apply \cref{dual1} and \cref{iden1}.
\end{proof}
\end{Coroll1}

We now show that, as in \cref{bimodstruc}, the prime ideal that annihilates the $\mathrm{Ext}$ bimodule on the left lies in the clique of $P$ when the clique of $P$ is localizable.

\begin{Lemma1} \label{dual3}
Let $R$ be an Auslander-Gorenstein grade-symmetric ring, and let $P$ be a prime ideal with a localizable clique $C$ with corresponding Ore set $X$. Suppose there exists a prime ideal $Q$ such that $\mathrm{Ext}^{j_{R}(R/P)}_{R}(R/P,R)$ is an $(R/Q,R/P)$-bimodule which is torsionfree on both sides. Then, $Q$ belongs to the clique of $P$.
\begin{proof}
By \cref{flatten}, it is true that $$R_{X}\otimes_{R} \mathrm{Ext}^{j_{R}(R/P)}_{R}(R/P,R)\cong \mathrm{Ext}^{j_{R}(R/P)}_{R_{X}}(R_{X}/P_{X}, R_{X}).$$ By \cref{locduality}, $R_{X}/P_{X}$ is a holonomic $R_{X}$-module. Thus, the module on the right hand side of the above formula is non-zero. This implies that $Q\cap X = \emptyset$, and thus implies that $Q$ is contained in some prime lying in the clique of $P$ (by the definition of a  localizable clique, the only maximal ideals of $R_{X}$ are those of the form $S_{X}$, as $S$ varies across $C$). As a left module, $\mathrm{Ext}^{j_{R}(R/P)}_{R}(R/P,R)$ is pure (\cref{extpure}) of grade $j_{R}(R/P)$. As it is also torsionfree as an $R/Q$-module, it follows that $j_{R}(R/Q)=j_{R}(R/P)$. Thus $Q$ itself belongs to the clique of $P$. 
\end{proof}
\end{Lemma1}

We summarize our results in the following theorem:

\begin{Theorems1} \label{bimodstruc2}
Let $A$ be a $k$-algebra with an Auslander, Cdim-symmetric dualizing complex of the form $R={^{\phi}A^{1}}[n]$ for some $\phi \in \mathrm{Aut(A)}$ which satisfies the trace property between classes of prime ideals $\mathcal{F}_{A}$ and $\mathcal{F}'_{A}$. If $P$ is a prime ideal in $\mathcal{F}_{A}$ with a localizable clique $C$ with corresponding Ore set $X$, then there exists a prime $D(P)\in C$ such that  $\mathrm{Ext}^{j_{A}(A/P)}_{A}(A/P,A)$ is an $(A/D(P),A/P)$-bimodule which is torsionfree and finitely generated on both sides. In this case, we have isomorphisms of $(A,A)$-bimodules: 

\begin{align*}
 \mathrm{Ext}^{j_{A}(A/P)}_{A}(A/P,A)\otimes_{A} A_{X} &\cong A_{X}\otimes_{A} \mathrm{Ext}^{j_{A}(A/P)}_{A}(A/P,A) \\
 & \cong \mathrm{Ext}^{j_{A}(A/P)}_{A_{X}}(A_{X}/P_{X},A_{X}). \\
\end{align*}

\begin{proof}
This follows from \cref{dual2} and \cref{dual3}. We can then argue as in \cref{bimod2}, or alternatively just use \cref{dual1}.
\end{proof}
\end{Theorems1}

The above result fulfills a promise made in the remark after \cref{bimodstruc}: for \cref{bimodstruc2}, we only need the second layer condition to be satisfied locally. 

\section{Localizing $\mathrm{Ext}^{i}_{R}(R/P,M)$}

We begin this section with two definitions.

\begin{Definitions1}
Let $R$ be a grade-symmetric Auslander-Gorenstein ring. A prime ideal $P$ is \textit{right correct} if $\mathrm{Ext}^{i}_{R}(R/P,R)$ is torsion as an $R/P$-module whenever $i\gneq j_{R}(R/P)$.
\end{Definitions1}

\begin{Definitions1} \label{defofloccor}
Let $R$ be a grade-symmetric Auslander-Gorenstein ring. A prime ideal $P$ is \textit{right localizably correct} if
\begin{itemize}
\item[1.)] The clique of $P$ is localizable, with corresponding Ore set $X$.
\item[2.)] $P$ is right correct.
\item[3.)] There are isomorphisms of $(R_{X},R_{X})$-bimodules: 
\begin{align*}
\mathrm{Ext}^{j_{R}(R/P)}_{R}(R/P,R)\otimes_{R} R_{X} & \cong R_{X}\otimes_{R} \mathrm{Ext}^{j_{R}(R/P)}_{R}(R/P,R) \\
& \cong   \mathrm{Ext}^{j_{R}(R/P)}_{R_{X}}(R_{X}/P_{X},R_{X}).\\
\end{align*}
\end{itemize}
A clique of prime ideals $C$ is \textit{right localizably correct} if all of the prime ideals in $C$ are right localizably correct.
\end{Definitions1}

There are obvious definitions of \textit{left correct} and \textit{left localizably correct} prime ideals, where we instead choose to work in the category of left modules. A prime ideal is \textit{correct} (resp.~\textit{localizably correct}) if it is both left and right correct (resp.~localizably correct). There is the obvious definition of a localizably correct clique.

\begin{Proposition1}
The following are examples of correct and localizably correct prime ideals:
 \begin{itemize}
\item [a.)] If $R$ is an Auslander-Gorenstein, grade-symmetric, weakly bifinite ring , then every prime ideal is correct. Furthermore, if $R$ satisfies the second layer condition, then every prime ideal with a localizable clique is localizably correct.
\item [b.)] Let $R$ be a $k$-algebra with an Auslander, Cdim-symmetric dualizing complex of the form $^{\phi}R^{1}[n]$ for some $\phi \in \mathrm{Aut(R)}$. Suppose $^{\phi}R^{1}[n]$ has the trace property between classes of prime ideals $\mathcal{F}_{R}$ and $\mathcal{F}'_{R}$. Then, every prime in $\mathcal{F}'_{R}$ (resp. $\mathcal{F}_{R}$) is right correct (resp. left correct).  If a prime $P$ has a localizable clique, then it is localizably correct.
\item [c.)] If $R$ is a grade-symmetric, Auslander-Gorenstein FBN ring, and $P$ is a prime ideal with a localizable clique $C$, then $P$ is localizably correct. 
\end{itemize}
\begin{proof}
The first part of a.) is precisely the content of \cref{onlyone}, b.) follows from \cite[Proposition 3.3]{JY}, and c.) is a consequence of \cref{localizeext}.

In cases $a.)$ and $b.)$, if $P$ has a localizable clique, then we can appeal to \cref{bimod2} and \cref{bimodstruc2} to see that it is localizably correct.
\end{proof}
\end{Proposition1}

\begin{Exam}
An interesting class of localizably correct prime ideals is provided by the work of \v{S}irola (\cite{SIR}). By \cite[Proposition 4.5 and Corollary 0.3]{JY}, it follows that $\mathcal{U}(\mathfrak{g})$, the universal enveloping algebra of a finite dimensional Lie algebra, satisfies the hypotheses in part $b.)$ of the above proposition. It is well known that every clique of prime ideals in the universal enveloping algebra of a solvable Lie algebra over $\mathbb{C}$ is classically localizable (\cite[Theorem A.3.5]{JAT}). \v{S}irola, in \cite[Corollary 1.5]{SIR}, gives a semisimple example: he shows that every height $1$ prime ideal in the universal enveloping algebra of $\mathfrak{sl_{2}}$ over $\mathbb{C}$ is localizable.
\end{Exam}

\begin{Ques}
If $P$ is a prime ideal with clique $C$, and $P$ is localizably correct, is every prime in the clique of $P$ localizably correct? 
\end{Ques}
We suspect that the answer to the above question is yes, but have no idea how to prove it. 

We will now use the structure of $\mathrm{Ext}^{j_{R}(R/P)}_{R}(R/P,R)$ to localize $\mathrm{Ext}^{i}_{R}(R/P,M)$; the spectral sequence in \cref{iss1} is the key to our arguments.

\begin{Lemma1} \label{iso1}
Let $R$ be a grade-symmetric Auslander-Gorenstein ring, and let $P$ be a correct prime ideal. Let $N$ be a right module of finite flat dimension. Then for all $i$, there is an isomorphism of right $R/P$-modules: $$\mathrm{Tor}_{i}^{R}(N,\mathrm{Ext}^{j_{R}(R/P)}_{R}(R/P,R)\otimes_{R} \mathcal{Q}(R/P))\cong \mathrm{Ext}^{j_{R}(R/P)-i}_{R}(R/P,N)\otimes_{R} \mathcal{Q}(R/P).$$
\begin{proof}
Recall, for suitable $M$ and $N$, the spectral sequence given in \cref{iss1}: $$E^{2}_{p,q}= \mathrm{Tor}_{p}^{R}(N, \mathrm{Ext}_{R}^{q}(M,R))\Rightarrow \mathrm{Ext}_{R}^{q-p}(M,N).$$  Taking $M=R/P$, and noting that the resulting modules are $R/P$-modules, \cref{flatspecseq} tells us that we can safely tensor with $\mathcal{Q}(R/P)$ to get a convergent spectral sequence  $$E^{2}_{p,q}= \mathrm{Tor}_{p}^{R}(N, \mathrm{Ext}_{R}^{q}(R/P,R)\otimes_{R} \mathcal{Q}(R/P))\cong \mathrm{Tor}_{p}^{R}(N, \mathrm{Ext}_{R}^{q}(R/P,R))\otimes_{R} \mathcal{Q}(R/P)$$ $$\Rightarrow \mathrm{Ext}_{R}^{q-p}(R/P,N)\otimes_{R} \mathcal{Q}(R/P).$$ The fact that $P$ is correct implies that $\mathrm{Ext}_{R}^{q}(R/P,R)\otimes_{R} \mathcal{Q}(R/P)\neq 0$ if and only if $q=j_{R}(R/P)$. The above spectral sequence thus collapses to give us the proposed isomorphism.
\end{proof}
\end{Lemma1}

The above lemma is the best we can do without any assumptions regarding the localizability of prime ideals. By \cref{dual1}, $\mathrm{Ext}^{j_{R}(R/P)}_{R}(R/P,R)\otimes_{R} \mathcal{Q}(R/P)$ is often isomorphic to a direct sum of copies of the unique simple module over $\mathcal{Q}(R/\phi(P))$ as a left module, where $\phi$ is as it is there. We are therefore interested in understanding $\mathrm{Tor}_{i}^{R}(N,\mathcal{Q}(R/\phi(P)))$. We have no idea how to do this without additional localizability assumptions.

\begin{Theorems1} \label{locmod}
Let $R$ be a grade-symmetric Auslander-Gorenstein ring, and let $P$ be a right localizably correct prime ideal with clique $C$ with corresponding Ore set $X$. Let $M$ be a module of finite flat dimension. Then, for all $i$ there is an isomorphism of right $R/P$-modules: $$\mathrm{Ext}^{i}_{R}(R/P,M)\otimes_{R} R_{X}\cong \mathrm{Ext}^{i}_{R_{X}}(R_{X}/P_{X},M_{X}).$$
\begin{proof}
Since $P$ is correct, \cref{iso1}  gives us the following isomorphism for $R$-modules $N$ of finite flat dimension: $$\mathrm{Tor}_{i}^{R}(M,\mathrm{Ext}^{j_{R}(R/P)}_{R}(R/P,R)\otimes_{R} \mathcal{Q}(R/P))\cong \mathrm{Ext}^{j_{R}(R/P)-i}_{R}(R/P,M)\otimes_{R} \mathcal{Q}(R/P).$$ Note that tensoring an $R/P$-module with $\mathcal{Q}(R/P)$ is the same thing as tensoring it with $R_{X}$.  Since $P$ is right localizably correct, we therefore get the following isomorphism: $$\mathrm{Tor}_{i}^{R}(M,\mathrm{Ext}^{j_{R}(R/P)}_{R_{X}}(R_{X}/P_{X},R_{X}))\cong \mathrm{Ext}^{j_{R}(R/P)-i}_{R}(R/P,M)\otimes_{R} R_{X}.$$ However, $$\mathrm{Tor}_{i}^{R}(M,\mathrm{Ext}^{j_{R}(R/P)}_{R_{X}}(R_{X}/P_{X},R_{X}))\cong \mathrm{Tor}_{i}^{R_{X}}(M_{X},\mathrm{Ext}^{j_{R}(R/P)}_{R_{X}}(R_{X}/P_{X},R_{X})).$$ Thus, the following isomorphism holds:  $$\mathrm{Tor}_{i}^{R_{X}}(M_{X},\mathrm{Ext}^{j_{R}(R/P)}_{R_{X}}(R_{X}/P_{X},R_{X}))\cong \mathrm{Ext}^{j_{R}(R/P)-i}_{R}(R/P,M)\otimes_{R} R_{X}.$$

However, we can also work in the category of $R_{X}$-modules. Using \cref{iss1} and \cref{locduality}, we see that there is an isomorphism $$\mathrm{Tor}_{i}^{R_{X}}(M_{X},\mathrm{Ext}^{j_{R}(R/P)}_{R_{X}}(R_{X}/P_{X},R_{X}))\cong \mathrm{Ext}^{j_{R}(R/P)-i}_{R_{X}}(R_{X}/P_{X},M_{X}).$$

Thus, we conclude that the proposed isomorphism holds.
\end{proof}
\end{Theorems1}

\section{Minimal Injective Resolutions}

Let $R$ be a noetherian ring, and let $M$ be a right $R$-module. If $P$ is a prime ideal of $R$, define the $i^{th}$ \textit{Bass} number of $M$ at $P$, $\mu^{i}(P,M)$, to be the multiplicity of $E_{P}$ in the $i^{th}$ term of the minimal injective resolution of $M$. Note that $\mu^{i}(P,M)$ is just a cardinal - it need not be finite.

The following well known result gives us information about $\mu^{i}(P,M)$. For a proof, we direct the reader towards \cite[Lemma 2.5]{JY}. 

\begin{Lemma1} \label{injnum}
 Let $R$ be a right noetherian ring, and $P$ a prime ideal in $R$. Let $M$ be an $R$-module. Then, $\mu^{i}(P,M)$ is the reduced rank of $\mathrm{Ext}^{i}_{R}(R/P,M)$ as an $R/P$-module. In particular, $\mu^{i}(P,M)\neq 0$ if and only if $\mathrm{Ext}_{R}^{i}(R/P,M)$ is not torsion as an $R/P$-module.
\end{Lemma1}

If $A$ is a commutative noetherian ring, then every indecomposable injective module is tame. Given a finitely generated $A$-module $M$, it is possible to describe the minimal injective resolution of $M$ in terms of geometric data arising from $M$. In particular, if $\mathrm{I}$ is the minimal injective resolution of $M$, the following facts are well known. Recall that a commutative noetherian ring is said to be \textit{Gorenstein} if it has finite injective dimension when viewed as a module over itself. If $A$ is a commutative noetherian local ring, and $M$ is a finitely generated $A$-module, the $i^{th}$ \textit{Betti} number of $M$, $\beta_{i}(M)$ is the rank of $P_{i}$, where $\mathrm{\mathbf{P}}:= ...\rightarrow P_{1}\rightarrow P_{0}\rightarrow 0$ is the minimal projective resolution of $M$.

\begin{itemize}
 \item Let $A$ be a commutative noetherian local ring. Let $P$ be a prime ideal in $A$. Then $\mu^{i}(P,M)\neq 0$ for some $i$ if and only if $P\in \supp{M}$ (\cite[Proposition 3.1.14]{BH}).
\end{itemize}

If $A$ is Gorenstein, and $M$ is assumed to have finite injective dimension, then more is true:

\begin{itemize}
\item If $P\in \supp{M}$, then $\mu^{i}(P,M)\neq 0$ if and only if $\depth{M_{P}}\leq i\leq \id{M_{P}}$ (\cite[Theorem 2.5]{FOX1}).
\item If $P\in \supp{M}$, then $\mu^{i}(P,M) = \beta_{\id{A/P}-i}(M_{P})$ (\cite[Corollary 3.2]{FOX1}).
\end{itemize}

We will attempt to generalize these facts to some families of noncommutative noetherian rings. Our main technique is going to be localization. Indeed, much of what we are going to prove will be an immediate consequence of what we have already done in this paper.

The reader should note one critical difference between what we are going to prove and what happens in the commutative case. For noncommutative rings, we are unable to determine when a particular prime $P$ appears in the minimal injective resolution of a module; the best we can do is determine precisely when \textit{some} prime from a given clique appears. Essentially this is because for noncommutative rings, a clique, as opposed to a single prime ideal, represents a single `blurred point' in the corresponding `noncommutative geometry'.

We first consider the case of FBN rings.

\begin{Theorems1} \label{resinFBN}
Let $R$ be a grade-symmetric Auslander-Gorenstein FBN ring, and let $C$ be a localizable clique with associated Ore set $X$. Let $M$ be a finitely generated right module of finite injective dimension, and let $\mathrm{I}$ be the minimal injective resolution of $M$. Then, $M_{X} \neq 0$ if and only if there exists some prime $Q\in C$ such that $\mu^{i}(Q,M)\neq 0$ for some $i$. Moreover, if $M_{X} \neq 0$, then for every $i$ between $\depth{M_{X}}$ and $\id{M_{X}}$, there exists a prime $Q\in C$ such that $\mu^{i}(Q,M)\neq 0$. 
\begin{proof}
Every simple left and right module of $R_{X}$ is a summand of $R_{X}/P_{X}$ for some prime $P\in C$, and by \cref{locduality}, $R_{X}$ is an AS-Gorenstein ring. Thus, we can apply \cref{niceclique}: given $i$ such that $\depth{M_{X}}\leq i\leq \id{M_{X}}$, there exists some prime $Q\in C$ such that $\mathrm{Ext}^{i}_{R_{X}}(R_{X}/Q_{X},M_{X})\neq 0$. 

By \cref{localizeext}, for all $i$ there is an isomorphism $$ \mathrm{Ext}^{i}_{R}(R/Q,M)\otimes_{R} R_{X}\cong \mathrm{Ext}^{i}_{R_{X}}(R_{X}/Q_{X},M_{X}).$$

Thus, for all $i$ between $\depth{M_{X}}$ and $\id{M_{X}}$, there exists a $Q\in C$ such that $\mathrm{Ext}^{i}_{R}(R/Q,M)$ is not torsion as an $R/Q$-module. The result now follows from \cref{injnum}.
\end{proof}
\end{Theorems1}

We now move on to more general rings. 

\begin{Theorems1} \label{resin1}
Let $R$ be a grade-symmetric Auslander-Gorenstein ring of dimension $d$. Let $M$ be a finitely generated module of finite injective dimension. Let $C$ be a right localizably correct clique with corresponding Ore set $X$. Then, $M_{X}\neq 0$ if and only if there exists a prime ideal $P\in C$ such that $\mu^{i}(P,M)\neq 0$ for some $i$.
\begin{proof}
We will again need to use \cref{locduality}. If there is a prime $P\in C$ such that $E_{P}$ appears as a summand of some term in the minimal injective resolution of $M$, then by \cref{injnum} and \cref{locmod}, $\mathrm{Ext}^{i}_{R_{X}}(R_{X}/P_{X},M_{X})\neq 0$ for some $i$. Thus $M_{X}\neq 0$. 

Conversely, if $M_{X}\neq 0$, then there exists a prime $Q\in C$ such that $M_{X}\otimes_{R_{X}} R_{X}/Q_{X}\neq 0$ --- this is immediate, since $M_{X}$ has a simple quotient. By \cref{locduality}, we see that $\mathrm{Ext}^{d}_{R_{X}}(R_{X}/D(Q)_{X},M_{X})\neq 0$. We can now use \cref{locmod} to complete our argument.
\end{proof}
\end{Theorems1}

A minor strengthening of our hypothesis allows us to recover the no-holes theorem as well.

\begin{Theorems1} \label{resin2}
Let $R$ be a grade-symmetric Auslander-Gorenstein ring. Let $M$ be a finitely generated module of finite injective dimension. Let $C$ be a right localizably correct clique with corresponding Ore set $X$. Suppose either of the following situations occur: 

\begin{itemize}
\item[a.)] $C$ is finite.
\item[b.)] $R$ contains an uncountable set $F$ of central units such that the difference of two distinct elements of $F$ is a unit, and the set of prime ideals $\{Q_{X}|Q\in C\}$ in $R_{X}$ satisfies the generic regularity condition.
\end{itemize}

If $M_{X}\neq 0$, then for every $i$ between $\depth{M_{X}}$ and $\id{M_{X}}$, there exists a prime ideal $Q\in C$ such that $\mu^{i}(Q,M)\neq 0$.
\begin{proof}
The proof is identical to the argument used in \cref{resinFBN}. We will need to use \cref{niceclique}, and \cref{locmod} instead of \cref{localizeext}. Note that cliques in noetherian rings are countable, by \cite[Theorem 16.23]{GW}.
\end{proof}
\end{Theorems1}

\begin{Exam}
In general, we cannot see inside the clique; determining when the tame indecomposable injective corresponding to \textit{some} prime in the clique appears as a summand of some term of a minimal injective resolution is the best that we can do through localization alone. Let us work through a particular example. Let $\mathfrak{g}$ be the unique non-abelian $2$-dimensional complex Lie algebra, with basis $x$ and $y$ and bracket $[x,y]=x$. For every $\alpha\in \mathbb{C}$, let $M_{\alpha}=\mathcal{U}(\mathfrak{g})/(x, y-\alpha)$. We will compute the minimal injective resolution of $M_{\alpha}$ over $\mathcal{U}(\mathfrak{g})$. Let $E_{\alpha}$ denote the injective hull of $M_{\alpha}$.

The clique of the prime ideal $(x, y-\alpha)$ is $\{(x,y-\alpha + n)\}_{n\in \mathbb{Z}}$. Thus the cliques of $\mathcal{U}(\mathfrak{g})$ are in bijection with $\mathbb{C}/\mathbb{Z}$. 

It is immediate that the only prime ideals whose corresponding indecomposable injectives appear in the minimal injective resolution of $M_{\alpha}$ are those in the clique of $(x, y-\alpha)$. Let us denote the Ore set corresponding to this clique by $X$. It is also true that only tame indecomposable injective modules appear in the minimal injective resolution of $M_{\alpha}$ (see remark after \cref{localizeext}). Thus we can completely determine the minimal injective resolution of $M_{\alpha}$ by working in $\mathcal{U}(\mathfrak{g})_{X}$.

We now make use of \cref{locduality}. The map $D$ is, in this case, well known: $D((x,y - \alpha +n))= (x, y - \alpha + n -1)$. If $I$ and $J$ are ideals in any ring $R$, $\mathrm{Tor}_{1}^{R}(R/I,R/J) = I\cap J/IJ$. Using this fact and \cref{locduality}, it easily follows that the minimal injective resolution of $M_{\alpha}$ is  $$0\rightarrow M_{\alpha}\rightarrow E_{\alpha}\rightarrow E_{\alpha}\oplus E_{\alpha + 1}\rightarrow E_{\alpha+1} \rightarrow 0.$$
\end{Exam}

\begin{Rem}
The above example suggests that the prime ideals which appear in the minimal injective resolution of a module are in some way determined by the internal link structure of cliques. This is a suggestive idea, but it is also clear that the injective dimension of a module plays some role in truncating the potential choices as well. We do not have a precise hypothesis about what is going on, but consider this as an interesting question for future research.
\end{Rem}

\subsection*{Multiplicities}

We now consider multiplicities. Recall that over a commutative Gorenstein ring, it is possible to precisely determine the multiplicity of $E_{P}$ in the minimal injective resolution of a finitely generated module. While, in general, we are unable to determine matters to that level of precision in the noncommutative case, we can at least conclude that the multiplicities are finite.

\begin{Proposition1} \label{finmult}
Let $R$ be a weakly bifinite grade-symmetric Auslander-Gorenstein ring. Let $M$ be a finitely generated module of finite injective dimension. Then, given a prime $P$, $\mu^{i}(P,M)$ is finite.
\begin{proof}
From the proof of \cref{iso1}, we see that $$\mathrm{Tor}_{i}^{R}(M,\mathrm{Ext}^{j_{R}(R/P)}_{R}(R/P,R)\otimes_{R/P} \mathcal{Q}(R/P))\cong \mathrm{Ext}^{j_{R}(R/P)-i}_{R}(R/P,M)\otimes_{R} \mathcal{Q}(R/P).$$ As $M$ is finitely generated and $R$ is weakly bifinite, the reduced rank of $\mathrm{Ext}^{i}_{R}(R/P,M)$ is finite for all $i$. \cref{injnum} completes the argument.
\end{proof}
\end{Proposition1}

\subsubsection*{Artin-Rees Rings}

There is a class of rings for which we can gain a better understanding of multiplicities. Recall that a noetherian ring is an Artin-Rees (AR) ring if every prime ideal satisfies the Artin-Rees property (see \cite[Chapter 13]{GW} for details on the Artin-Rees property in noncommutative ring theory). The class of noetherian AR rings in which every prime ideal is completely prime has good geometric properties (see, for example, \cite{SIR}). 

If $R$ is a scalar-local noetherian ring with unique maximal ideal $\mathfrak{m}$, then every finitely generated  $R$-module has a minimal free resolution --- this is true because $R$ is semi-perfect, and over a scalar-local ring every finitely generated projective module is free. We extend the definition of a Betti number to this noncommutative setting. If $...\rightarrow F_{i}\rightarrow ... \rightarrow F_{1}\rightarrow F_{0}\rightarrow 0$ is the minimal free resolution of a module $M$, we define $\beta_{i}(M)$, the $i^{th}$ \textit{Betti} number of $M$, to be the rank of $F_{i}$. If $k:=R/\mathfrak{m}$, then $\beta_{i}(M) = \text{dim}_{k}\mathrm{Tor}^{R}_{i}(k,M).$

An artinian ring $B$ is said to be \textit{weakly symmetric} (\cite[Definition 1.6]{QZ2}) if given rings $L$ and $C$, and bimodules $_{C}M_{B}$ and $_{L}N_{B}$ (resp.~$_{B}M_{C}$ and $_{B}N_{L}$) of finite length on both sides, $\mathrm{Hom}_{B}(M,N)$ is an $(L,C)$ (resp.~ $(C,L)$) bimodule of finite length on both sides.

\begin{Theorems1} \label{algar}
Let $R$ be a complete scalar-local Auslander-Gorenstein algebra of dimension $d$ over a field $k$ with maximal ideal $\mathfrak{m}$, and suppose that every prime ideal is completely prime and satisfies the Artin-Rees property. Suppose further that $R/\mathfrak{m}$ is weakly symmetric. Then, given a finitely generated module $M$ of finite injective dimension, $$\mu^{i}(P,M)=\beta_{j_{R}(R/P)-i}(M_{P}).$$
\begin{proof}
Since $R$ is scalar-local and Auslander-Gorenstein, it is AS-Gorenstein. It follows from \cite[Theorem 0.1]{CWZ} that $R$ is grade-symmetric, weakly bifinite, and has finite cohomological dimension. Since $R$ is an $AR$ ring, it satisfies the second layer condition, and every prime ideal is classically localizable; in particular, every  clique in $R$ is a singleton. Thus, \cref{resin2} applies. 

Let $P$ be a prime ideal in $R$, and let $X$ be $R\  \backslash \  P$ (this is the Ore set corresponding to the clique of $P$ in this case). Using \cite[Proposition 4.5]{JY}, we see that our hypotheses imply that $R$ has the trace property for prime ideals. Since every clique in $R$ is a singleton, it follows from the trace property and \cref{dual3} that $\mathrm{RHom}_{R}(R/P,R)$ is a dualizing complex between $R/P$ and itself; in particular $\mathrm{Ext}^{j_{R}(R/P)}_{R}(R/P,R)$ is annihilated by $P$ and is finitely generated and torsionfree as an $R/P$-module on both sides. 

Since every prime ideal in $R$ is completely prime, $R_{X}/P_{X}$ is a division ring.  \cite[Proposition 3.3]{JY} then implies that $$B:= R_{X}\otimes_{R} \mathrm{Ext}^{j_{R}(R/P)}_{R}(R/P,R)\otimes_{R} R_{X}\cong \mathrm{Ext}^{j_{R}(R/P)}_{R_{X}}(R_{X}/P_{X},R_{X})$$ is an $(R_{X}/P_{X},R_{X}/P_{X})$-bimodule of length $1$ on both sides.

An easy computation now tells us that $\mathrm{Tor}^{R_{X}}_{i}(M_{X},B)$ has dimension $\beta_{i}(M_{X})$ --- to see this, take a minimal projective resolution of $M_{X}$ over $R_{X}$, and use the fact that it is \textit{minimal}. 

However, by the collapsing of the Ischebeck spectral sequence in \cref{iss1} (see the proof of \cref{locmod} for details), $\mathrm{Tor}^{R_{X}}_{i}(M_{X},B)\cong \mathrm{Ext}^{j_{R}(R/P)-i}_{R_{X}}(R_{X}/P_{X},M_{X})$. It follows that $\mu^{i}(P,M)=\beta_{j_{R}(R/P)-i}(M_{P}).$
\end{proof}
\end{Theorems1}

\begin{Rem}
Completed group algebras of nilpotent uniform pro-p groups over $\mathbb{F}_{p}$ satisfy the hypotheses of the following theorem; this follows from the work of Ardakov (\cite{ARD1}). 

In fact, the above theorem continues to hold if we only assume that $R$ is an algebra over a principal ideal domain; in the forthcoming paper \cite{VYAS2}, we observe that with minor modifications in proof, the results involving dualizing complexes used in the proof of \cref{algar} continue to hold in this greater generality.
\end{Rem}

\section{An Example}

The reader may wonder whether the homological conditions we imposed on our rings in \textsection{3} are necessary. In this section, we show that we cannot proceed naively; some additional conditions are needed to gain traction over our rings, even in superficially nice situations. 

\begin{Proposition1}
There is a noetherian $\mathbb{C}$-algebra $R$ such that
\begin{itemize}
 \item[a.)] $R$ has Krull dimension $1$,
 \item[b.)] $R$ has classical Krull dimension $1$,
 \item[c.)] $R$ has global dimension $1$,
 \item[d.)] $R$ satisfies the second layer condition,
 \item[e.)] Every clique in $R$ is classically localizable.
\end{itemize}
 There is a prime ideal $Q$ of $R$, with clique $C$ with corresponding Ore set $X$, and a finitely generated uniform $R$-module $M$ such that $\mathrm{Ext}^{1}_{R}(R/Q,M)\otimes_{R} R_{X}\neq 0$, but $\mathrm{Ext}^{1}_{R_{X}}(R_{X}/Q_{X},M_{X})= 0$.
\begin{proof}
Our example is the one constructed by Goodearl and Schofield in \cite{GS}.

In the example after \cite[Corollary 4]{GS}, Goodearl and Schofield construct a right and left noetherian ring $R$ with right and left Krull dimension $1$ such that there is a simple module with a non-artinian finitely generated essential extension. Let us show that we choose $R$ in such a way that $R$ satisfies the second layer condition and every clique in $R$ is classically localizable.

The rings constructed in \cite{GS} are triangular matrix rings of the form $R= \bigl(\begin{smallmatrix}
S& T \\ 0& T
\end{smallmatrix} \bigr)$, where $S$ and $T$ are noetherian rights of Krull dimension $1$ such that $S$ is a subring of $T$ with $T$ finitely generated as a left $S$-module. In particular, for the examples constructed in \cite{GS}, $S\cong F[X]$ and $T\cong E[X]$ for a suitable extension $F\subset E$ of division rings. Studying \cite[Corollary 4]{GS}, we see that we can assume that $F$ contains $\mathbb{C}$; in particular, we can assume that $R$ is an algebra over an uncountable field.

Let $N$ denote the nilradical of $R$. It follows that $N^{2}=0$, and $R/N\cong F[X]\times E[X]$. It is now immediate that $R$ has left and right Krull dimension $1$. By \cite[Proposition 16.9]{GW}, $\spec{R}$ satisfies the generic regularity condition. It is clear that the classical Krull dimension of $R$ is $1$.

Let us show that $R$ has global dimension $1$. By \cite[Proposition 5.1]{MR}, we see that the global dimension of $R$ is less than or equal to $2$; further, by \cite[Corollary 1.14]{MR}, the global dimension of $R$ is the supremum of the projective dimensions of its simple modules. Let $P= \bigl(\begin{smallmatrix} 0& T \\ 0& T \end{smallmatrix} \bigr)$, and $Q= \bigl(\begin{smallmatrix}
S& T \\ 0& 0
\end{smallmatrix} \bigr)$. It is not to hard to see that $P$ and $Q$ are the only two minimal prime ideals of $R$, and $PQ=0$. Every simple module must therefore have annihilator containing $P$ or $Q$. If $D$ is a division ring, the standard division algorithm immediately implies that $D[X]$ is a principal right and left ideal ring. Thus, every simple right $R$-module is either isomorphic to $$ \bigl(\begin{smallmatrix}
S& T \\ 0& T
\end{smallmatrix} \bigr) / \bigl(\begin{smallmatrix}
S& T \\ 0& fT
\end{smallmatrix} \bigr),$$ or $$ \bigl(\begin{smallmatrix}
S& T \\ 0& T
\end{smallmatrix} \bigr) / \bigl(\begin{smallmatrix}
gS& T \\ 0& T
\end{smallmatrix} \bigr),$$ for some $f\in T$ or $g\in S$.

Since both $\bigl(\begin{smallmatrix}
S& T \\ 0& fT
\end{smallmatrix} \bigr)$ and  $\bigl(\begin{smallmatrix}
gS& T \\ 0& T
\end{smallmatrix} \bigr)$ are isomorphic to $R$ as right $R$-modules, it is immediate that the projective dimension of every simple $R$-module is less than or equal to $1$. Since $R$ is clearly not semisimple, it has global dimension equal to $1$.

We now show that $R$ satisfies the second layer condition. By \cite[Theorem 12.6]{GW}, it is enough to show that $D[X]$ satisfies the condition, where $D$ is a division ring. By \cite[Proposition 17.1]{GW}, every ideal in $D[X]$ contains a central element. Thus $D[X]$ has normal separation, and by \cite[Theorem 12.17]{GW}, we conclude that $D[X]$ satisfies the second layer condition. 

By \cite[Lemma 6]{W}, \cite[Theorem 7.1.5]{JAT}, \cite[Theorem 16.22]{GW}, and \cite[Corollary 14.6]{GW}, every clique in $R$ is classically localizable.

We now produce our counterexample. If $U$ is a uniform module over a right noetherian ring, it has a unique associated prime, say $P$. Define $\operatorname{pt}(U):=\{x\in U| xP=0\}$. Since there is a simple module with a non-artinian finitely generated essential extension, \cite[Theorem 3.5]{JAT2} and \cite[Corollary 3.6]{JAT2} tell us that there is a finitely generated, uniform $R$-module $M$ with unique associated prime $P$ such that
\begin{itemize}
 \item[i.)] $\operatorname{pt}(M)$ is torsion as a $R/P$-module.
 \item[ii.)] $M/\operatorname{pt}(M)$ is uniform with associated prime $Q$, and $MQ\leq \operatorname{pt}(M)$.
 \item[iii.)] $M/\operatorname{pt}(M)$ is torsionfree as an $R/Q$-module.
\end{itemize}

By \cite[Theorem 12.1]{GW}, we see that $Q\rightsquigarrow P$. Corresponding to the short exact sequence $$0\rightarrow \operatorname{pt}(M)\rightarrow M\rightarrow M/\operatorname{pt}(M)\rightarrow 0,$$ there is a long exact sequence $$0\rightarrow \mathrm{Hom}_{R}(R/Q,\operatorname{pt}(M))\rightarrow \mathrm{Hom}_{R}(R/Q,M)\rightarrow \mathrm{Hom}_{R}(R/Q,M/\operatorname{pt}(M))\rightarrow \mathrm{Ext}^{1}_{R}(R/Q, \operatorname{pt}(M)).$$ 

Since $Q$ is not contained in $P$, it follows that $\mathrm{Hom}_{R}(R/Q,\operatorname{pt}(M))\cong \mathrm{Hom}_{R}(R/Q,M)=0$. We see that $M/\operatorname{pt}(M)$ embeds as a submodule of $\mathrm{Ext}^{1}_{R}(R/Q,\operatorname{pt}(M))$. In particular, $$\mathrm{Ext}^{1}_{R}(R/Q,\operatorname{pt}(M))\otimes_{R} R_{X}\neq 0.$$ However, since $\operatorname{pt}(M)$ is torsion as a $R/P$-module, $\mathrm{Ext}^{1}_{R_{X}}(R_{X}/Q_{X},\operatorname{pt}(M)_{X})=0$.
\end{proof}
\end{Proposition1}

\begin{Rem}
In a forthcoming paper, \cite{VYAS1}, we use some of the ideas in the above proposition to give a homological reformulation of the link condition.
\end{Rem}

\section*{Acknowledgments}

This work was done as part of the author's PhD thesis, which was supervised by Simon Wadsley. The author would like to thank Simon Wadsley for his many helpful comments and suggestions while reading through many preliminary versions of this paper. He would like to thank James Zhang for making an observation that motivated a large part of \textsection 3. The author would also like to thank Jonathan Nelson for his careful proof-reading of a preliminary version of this paper, and his suggestions and helpful comments. The author thanks the Cambridge Commonwealth Trust and Wolfson College, Cambridge for their support over the years. The author would also like to thank the anonymous referee for their comments and suggestions.

\end{document}